\documentclass[10pt,a4paper]{article}
\usepackage{latexsym,amssymb}
\usepackage{enumerate, verbatim}
\newtheorem{theorem}{Theorem}
\newtheorem{lemma}{Lemma}
\newtheorem{prop}{Proposition}
\newtheorem{corollary}{Corollary}

\newcommand{\pf}{{\it Proof:\quad}}
\newcommand{\dne}{\hfill $\Box$\vspace{0.3cm}}

\newcommand{\Sn}{{\rm Sym}_{n}}

\newcommand{\gG}{\Gamma}

\newcommand{\ga}{{\alpha}}
\newcommand{\gb}{{\beta}}

\newcommand{\spa}{\!\!\!\!\!\phantom{|^{|^{|}}_{|_{|^{|}}}}\!\!\!\!\!}
\begin{document}

\title{Error Graphs and the Reconstruction \\ of Elements in Groups}
\author{
Vladimir I. Levenshtein\thanks{This research was supported by
the Russian Foundation for Basic Research (Grant 04-01-00112). }
%\and 
\\Keldysh Institute of Applied Mathematics,\\
Russian Academy of Sciences, Moscow, Russia\\
leven@keldysh.ru\\\\
\and
Johannes Siemons 
\\ School of Mathematics, \\University of East Anglia, Norwich, UK\\
j.siemons@uea.ac.uk}
\date{\scriptsize Accepted Version, November 20, 2008; printed \today}
\maketitle
\begin{abstract}

Packing and covering problems for metric spaces, and graphs in particular,  are of essential interest in combinatorics and coding theory.
They are formulated in terms of metric balls of vertices. We consider
a new problem in graph theory which is also based on the consideration of
metric balls of vertices, but which is distinct from the traditional packing and
covering problems. This problem is motivated by
applications in information transmission when redundancy of messages
is not sufficient for their exact reconstruction, and applications in computational
biology when one wishes to restore an evolutionary process. It can be
defined as the  reconstruction, or identification, of an unknown vertex in a given graph from a minimal number of vertices (erroneous or distorted patterns) in a metric ball of a given radius $r$ around the unknown vertex. For  this problem it is required to find minimum restrictions for such a reconstruction to be  possible and
also to find efficient  reconstruction algorithms under such minimal restrictions.

\medskip
In this paper we define error graphs and investigate  their basic properties. A particular class of error graphs occurs when the vertices of the graph are the elements of a group, and when the path metric is determined by a suitable set of group elements. These are the  undirected  Cayley graphs. Of particular interest is the transposition Cayley graph on the symmetric group which occurs in connection with the analysis of transpositional mutations in molecular biology~\cite{pev00,sa02}. We obtain a complete solution of the above problems for the transposition Cayley graph on the symmetric group.\newline

\phantom{XXXXXXXsX}Stirling Numbers


\end{abstract}

\section{Introduction: A Graph-Theoretical Approach to \\ Efficient
Reconstruction}

%{\scriptsize J: In this section there is little change,  mostly
%polishing of notation. I have mentioned 'error and distortion'
%already at the beginning so that the connection to efficient
%reconstruction is clearer.}

The problem of the efficient reconstruction of sequences was
introduced  in \cite{lev97,eros, eros2} as a problem in coding theory, and similar questions about the 
efficient reconstruction of integer partitions were considered in \cite{masi,presi}. 
In this paper we discuss a graph-theoretical setting in which efficient reconstruction problems can be studied as a uniform  theory.  

Let $\Gamma =(V, E)$
be a simple, undirected and  connected graph with vertex set $V$ and edge set $E$.
We regard the vertices in  $V$ as   units of information  in the given reconstruction problem, and for two vertices $x\neq y$ in $V$ we regard  $\{x,\,y\}$ as  an edge of $\gG$ if $y$ is obtained from $x,$ or vice versa $x$ from $y$, by a {\it single error}  or  {\it single distortion} of information. We might say that $x$ and $y$ are erroneous single error representations of each other, and that $\gG$ is a single error graph. The precise definitions can be found in Section 2.  The task of the reconstruction problem now  is to {\it restore}\, or {\it reconstruct}\, the original unit of information from sufficiently many erroneous representations of it. In other words, an unknown vertex $x$ in $\gG$ is to be identified by suitable knowledge about its neighbouring vertices in $\gG.$

We denote  the path distance between two vertices $x$ and $y$ of $\gG$  by $d(x,y)$
and we let %$S_r(x)=\{\,y\in V\,\,:\,\,  d(x,y)= r\,\} \mbox{\quad and\quad} 
$B_r(x)=\{\,y\in V\,\,:\,\,  d(x,y)\leq r\,\}$ be the
%{\it sphere }  and 
 ball  of radius $r$ centered at  $x.$
For given $r\geq 1$ denote by $N(\Gamma,r)$ the largest number $N$
such that there exist a set $A\subseteq V$ of size $N$ and two
vertices $x\neq y$ with $A\subseteq B_r(x)$ and $A\subseteq B_r(y).$
Thus any $N+1$ distinct vertices are contained in $B_r(x)$  for at
most one vertex $x$ while there are some $N$ vertices simultaneously contained in $B_{r}(x)$ and $B_{r}(y)$ for some $x\neq y.$ This means that an unknown vertex of $\gG$ can be identified, or 
reconstructed uniquely, by any set of $N(\Gamma,r)+1$ or more
distinct vertices at distance at most $r$ from the vertex, provides that such a set exists.

In graph theoretical terms we are therefore required, for an arbitrary  graph $\Gamma$ and an integer $r\geq 1,$  to determine the number 
\begin{equation}
\label{e001} N(\Gamma,r)=\max_{x,y\in V,\ x\neq y}| B_r(x)\cap
B_r(y)|
\end{equation}
and to construct an efficient  algorithm by which any  unknown
vertex $x$ in $V$ can be identified uniquely from an arbitrary  set of $N(\Gamma,r)+1$ vertices at distance $r$ or less from $x.$ Evidently we can assume that $r$ is at most $d(\Gamma), $ the diameter of $\gG.$ Throughout the paper we
assume that  $d(\Gamma)\ge 2$ and in particular $|V|\ge 3$. 

\bigskip
Problems of this kind have been solved for some graphs  and metric spaces of
interest in coding theory, and to give an impression of such results we review the example of Hamming spaces and Johnson spaces.  The Hamming space $F_q^n$ consists of
$q^n$ vectors of length $n$ over the alphabet $\{0,1,...,q-1\}$ with
metric $d(x,y)$ given by the number of coordinates in which the
vectors $x$ and $y$ differ. This metric space can be represented by
a graph $\gG$ whose vertices are the vectors of $F_q^n$ with two
vectors connected by an edge if and only if they differ in a single
coordinate. The path distance between two vertices then is the
Hamming distance between the  corresponding vectors. Therefore we can identify $F_q^n$ with this graph $\Gamma$. 
In \cite{lev97,eros, eros2}
it was shown that for any $n,$ $q$ and $r$ we have 
\begin{equation}
\label{e002} N(F_q^n,r)=q\sum_{i=0}^{r-1}\left({n-1} \atop{i}\right)
(q-1)^i\,\,. 
\end{equation}
Furthermore, any $x\in F_q^n$ can be reconstructed from $N=N(F_q^n,r)+1$
vectors of $B_r(x),$ written as the columns of a matrix, by applying the
majority algorithm to the rows of the matrix.

For any \,$1\le w \le n-1$\, the Johnson space $J_w^n$ consists of the ${n\choose
w}$ binary vectors in $F_2^n$ of length $n$ and Hamming weight $w$, where
distance is equal to half  the (even) Hamming distance in $F_2^n.$ This distance coincides with the minimal number
of coordinate transpositions needed to  transform  one vector into the
other. The Johnson space then can be viewed as a graph $\Gamma$ whose
vertices are the  vectors of $J_w^n$ with two vectors connected by
an edge if and only if one is obtained from the other by a
transposition of two coordinates. The path distance between two vertices
of $\Gamma$ then is the Johnson distance between the corresponding
vectors. Therefore we can identify $J_w^n$ with this graph $\Gamma$.
In \cite{lev97,eros} it was also shown that for any $n,$ $w$ and
$r$ we have \begin{equation}
\label{e003} N(J_w^n,r)=n\sum_{i=0}^{r-1}\left({w-1} \atop{i}\right)
\left ({n-w-1} \atop{i}\right) \frac 1{i+1}\,\,.
\end{equation}
Furthermore, any $x\in J_w^n$ can be reconstructed from
$N=N(J_w^n,r)+1$ vectors of $B_r(x),$ written as the columns of a matrix,
by applying a threshold algorithm to the rows of the matrix.

\bigskip
In the first part of this paper we make the notion of  error graphs precise  and  develop the theory needed to estimate $N(\gG,r)$ in  some general situations. In this respect our main results are  Theorems~\ref{th001}\, and~\ref{th002}\, which give  lower bounds for $N(\gG,1)$ and $N(\gG,2)$ in terms of other graph parameters. It may be useful to mention  that the idea of reconstructing a vertex in a given graph has nothing to do, a priory,  with the classical  Ulam problem of reconstructing a graph from the isomorphism classes of its vertex-deleted subgraphs. So we do not refer to the well-known and unresolved vertex-reconstruction problem.
Nevertheless, error graphs are such a general tool that even this problem can be phrased suitably a problem on error graphs.

In the second part of the paper we deal with error graphs for which the vertex set consists of  the elements of a group, and where the errors are defined  by a certain set of  group elements. Such graphs turn out to be  undirected Cayley graphs, and in Sections 4 and 5\,  we show that many important error graphs occur as Cayley graphs.  In Section 5 we discuss how transpositional errors in biological nucleotide sequences can be described as  errors in the transposition Cayley graph $\Sn(T)$ on the symmetric group. The remainder of the paper deals with this graph in particular. 

In Theorem~\ref{con1}\, we determine the full automorphism group of  the transposition Cayley graph $\Sn(T).$ 
The explicit value of $N(\Sn(T),r)$ can be found in Theorems~\ref{lem0008},\,\,\ref{th005}\, and \ref{th006}\, for $1\leq r\leq 3.$  To state the main result on $N(\Sn(T),r)$ for arbitrary $r\geq 1$ let $c(n,n-r)$ be the number of permutations on $\{1..n\}$ having exactly $n-r$ cycles. Thus the $c(n,n-r)$ are the signless Stirling numbers of the first kind. We also need the following {\it restricted Stirling numbers}: Let $c_{3^{1}}(n,n-r)$ be the number of permutations $g$ on $\{1..n\}$ having exactly $n-r$ cycles such that $1,\,2$ and $3$ belong to the same cycle of $g.$ The main result on $N(\Sn(T),r)$ is Theorem~\ref{th0009}. It shows that for all $r\geq 1$   %we have  
\begin{eqnarray}
\label{lastt} 
N(\Sn(T),r)&=&\sum_{i=0}^{r-1}\,\,c(n,n-i)  \,\nonumber\\
&+&  c_{3^{1}}(n,n-r) + c_{3^{1}}(n,n-(r+1))\,.\,\,\,\spa
\end{eqnarray}
for all  sufficiently large $n.$ Furthermore, the maximum $N(\Sn(T),r)=|B_{r}(x)\cap B_{r}(y)|$ occurs for any $x\neq y$ for which $x^{-1}y$ is a $3$-cycle on $\{1..n\}.$ We mention the connection between  this  theorem and  the  Poincar\'{e} polynomial  of $\Sn(T).$ When $\Gamma$ is an arbitrary finite graph and $v$ a vertex of $\gG$ let $c_{i}$ denote the number of vertices at distance $i$ from $v.$ Then 
$$\Pi_{\Gamma\!,v}\,(t):= \sum_{0\leq i}\, c_{i}t^{i}$$ is the {\it Poincar\'{e} polynomial } of $\Gamma$ at $v.$ When this polynomial is independent of $v$ we write simply 
$\Pi_{\Gamma}(t).$ For  the transposition Cayley graph $\Sn(T)$ the Poincar\'{e} polynomial is 
\begin{equation}\label{poin}
\Pi_{\,\Sn(T)}\,(t)=\sum_{i=0}^{n-1}\,\,c(n,n-i)\,t^{i}
\end{equation}

where the $c(n,n-i)$ are the Stirling numbers appearing in~~(\ref{lastt}). This shows that  the reconstruction parameters $N(\gG,r)$ are  related to important graph  invariants. 

In this paper we have avoided technical terminology as far as possible in order to make this material accessible to non-specialists. For the same reasons we have added a few key references to texts in computing and computational biology.

\bigskip
\section{Errors in Graphs}
We will now fix the notation used for the remainder. Let $\gG=(V,E)$ be a finite graph with vertex set $V$ and edge set $E.$ All edges are undirected and there are no multiple edges or loops.  Let $x,\,y$ be vertices. Then $x$ and $y$ are  {\it adjacent } to each other if $\{x,\,y\}$ is an edge. Further,  $d(x,y)$ denotes  the usual graph distance between the vertices, that is the length of a shortest path from $x$ to $y.$ Put $d(x,\,y)=\infty$ if $x$ and $y$ are in different components. For $i\geq 0$ we let 
$B_{i}(x):=\{ y\in V\,\,:\,\, d(x,y)\leq i\}$ and $S_{i}(x):=\{ y\in V\,\,:\,\, d(x,y)= i\}$ be the {\it ball } and {\it sphere} of radius $i$ around $x,$ respectively. 

We put $k_i(x)=|S_i(x)|$ and for  $y\in S_i(x)$ we set
\begin{eqnarray}\label{e00001}
c_i(x,y)&:=&|\{z\in S_{i-1}(x)\,\,:\,\, d(z,y)=1\}|\,,\spa\nonumber \\
a_i(x,y)&:=&|\{z\in S_i(x)\,\,: \,\,d(z,y)=1\}|\,,\spa\nonumber \\
b_i(x,y)&:=&|\{z\in S_{i+1}(x)\,\,:\,\, d(z,y)=1\}|\,.
\end{eqnarray}

It is clear that $b_0(x,y)=k_1(x),$ that $a_1(x,y)=a_{1}(y,x)$ is the number of
triangles over the vertices $x$ and $y$, and that $c_2(x,y)$ is the number
of common neighbours of $x$ and $y\in S_2(x).$ Let
\begin{eqnarray}
\label{e004} \lambda&=&\lambda(\Gamma)=\max_{x,y\in V,\, d(x,y)=1\,\,}a_1(x,y)\quad \spa \nonumber \\  \mu &=&\mu(\Gamma)=\max_{x,y\in V,\,d(x,y)=2\,\,}c_2(x,y).
\end{eqnarray}
Since $|B_r(x)\cap B_r(y)|>0$ for $x\neq y$ only if $
d(x,y)\le 2r$ we have
\begin{equation}
\label{e005} N(\Gamma,r)=\max_{1\le s \le 2r}N_s(\Gamma,r)
\end{equation}
where
\begin{equation}
\label{e006} N_s(\Gamma,r)=\max_{x,y\in V,\, d(x,y)=s}| B_r(x)\cap
B_r(y)|.
\end{equation}

In  particular, $N_1(\Gamma,1)=\lambda +2\,$ and $\,N_2(\Gamma,1)=\mu$ so that
\begin{equation}
\label{e007} N(\Gamma,1)=\max(\lambda +2, \mu).
\end{equation}
Finding or  estimating the value $N(\Gamma,r)$ for graphs of
interest in applications is  the main aim of our investigation here. We note the following general bounds for $N(\Gamma,r).$

\begin{lemma} \label{lem0001} \, Suppose that $x\neq y$ are vertices in the connected graph $\gG=(V,E)$ at distance $s=d(x,y)$ from each other. Let $r\geq 0$ be an integer.

\vspace{-0.3cm} (i) \,\, If $r\geq s$ then
$B_{r}(x)\cap B_{r}(y)=B_{r-s}(x)\,\cup\, [(B_{r}(x)\setminus B_{r-s}(x))\cap B_{r}(y)].\spa$ In particular, we have $N_{s}(\gG, r)\geq |B_{r-s}(x)|$ and 
\begin{equation}
\label{e0001}    N(\gG, r)\geq |B_{r-1}(x)| 
\end{equation}
for some $x\in V.$

(ii) \,\, If $r< s$ then
$B_{r}(x)\cap B_{r}(y)=B_{r}(y)\,\cap\, [B_{r}(x)\setminus B_{s-r-1}(x)]$.
\end{lemma}

\pf One should think of $B_{r}(x)\setminus B_{r-s}(x)$ as an annulus around $x.$ (i)~Starting on  a path of length $s$ from $y$ to $x$ any vertex in $B_{r-s}(x)$ can be reached by a further path of length at most $r-s$. The other statements are immediate from this.  (ii)~This is an direct  consequence of the triangle inequality.
\dne 

We set
\begin{equation}
\label{e008} k_i(\Gamma)=\max_{x\in V}k_i(x)\quad\mbox {and}\quad
k(\Gamma)=k_1( \Gamma).
\end{equation}
Then $\Gamma$ is {\it regular}\, of valency $k$ (or {\rm $k${\it -regular)\,}
if all its vertices have constant valency $k=k(\Gamma)$. A
$k$-regular graph is {\it distance-regular } if the numbers
$c_i(x,y)$ and $b_i(x,y)$ (and hence $a_i(x,y)=k-c_i(x,y)-b_i(x,y)$)
do not depend on $x\in V$ and $y\in S_i(x)$, for all
$i=0,1,...,d(\Gamma)$. A distance-regular graph of diameter $2$\, is
{\it strongly regular.} A good reference to strongly regular graphs is Chapter~21 \,in~\cite{Wilson} or also \cite{bcn}. In such a graph there are integers $\lambda$ and $\mu$ so that  any pair of vertices $x\neq y$ is simultaneously adjacent to exactly $\lambda$ vertices if $\{x,\,y\}$ is an edge, and to exactly $\mu$ vertices if $\{x,\,y\}$ is not an edge. Our use in \,(\ref{e004})\, of the symbols $\lambda$ and $\mu$ is therefore a natural extension to graphs which are not strongly regular.

Let ${\rm Aut}(\Gamma)$ be the automorphism group of $\Gamma.$ If $\gG$ is
vertex-transitive (that is, for any two vertices in $V$ there is an automorphism of $\Gamma$ mapping one onto the other) then
$k_i(x)=k_i(\Gamma)$ is constant for  all $x\in V$ and $i$. In
particular, such a graph is regular. However, even for vertex-transitive graphs
the $c_i(x,y)$ and $b_i(x,y)$ usually depend on $y\in S_i(x),$ and this
can cause difficulties in finding  $N(\Gamma,r).$ This phenomenon can be observed already on relatively small graphs, see the Remark following Lemma~\ref{lem0007}.

The Hamming and Johnson graphs are examples of error graphs  in which 
two vertices $x\neq y$ are joined by an edge if and only if there exists a
single error (the substitution of  a symbol or the transposition of two
coordinates, respectively) which transforms $x$ to $y$ and there exists
a single error which transform $y$ to $x.$ This observation leads to  a natural general theory of  single errors which we began in~\cite{eros}.  
For this we let $V$  be  a finite (or countable) set. A {\it single error } on $V$ is
an injection $h:V_{h}\rightarrow V$ defined on a non-empty subset
$V_{h}\subseteq V$ so that $h(x)\neq x$ for all $x\in V_{h}$. A
non-empty set $H$ of single errors will be called a {\it single
error set,}\, or just {\it error set, }
provided the following two properties hold:

\vspace{-0.5cm}
\begin{enumerate}[(i)]

\item For each $h\in H$ and $x\in V_{h}$ there exists some $g\in H$ so that
$h(x)\in V_{g}$ and $g(h(x))=x$, and
\item For all distinct pairs $x, y\in V$ there exist $x=x_1,x_2,...,x_m=y\in V$ and
$h_1,h_2,...,h_{m-1}\in H$ such that $x_{i+1}=h_i(x_i),$ for
$i=1,...,m-1$.
\end{enumerate}

%{\scriptsize $E:$ (ii) For all distinct pairs $x, y\in V$ there exist a path$x=x_1,x_2,...,x_m=y\in V$ and $h_1,h_2,...,h_{m-1}\in H$ such that $x_{i+1}=h_i(x_i),$ for $i=1,...,m-1$.\quad $J:$ We will mention the word path later, at the moment the edges aren't defined yet, and there is no graph as yet. Note also that I have simplified notation a little: once it is clear that $H$ always means  {\it single} errors we don't have to repeat it all the time. Again, English is a little more flexible than Russian (or German) would be.}  

For such a set $H$ we construct the {\it error graph} $\Gamma_H=(V, E)$ where  $E=\{\,\{x,h(x)\}\,\,:\,\,x\in V\mbox{\,\,and \,\,}
h\in H\,\}$. Note, by the conditions on $H$ we see that $\gG_{H}$ has no loops and 
that all edges are undirected. The condition (ii) says that there is a path between any two vertices,  and hence that $\gG_{H}$ is connected.
Furthermore, the usual path distance $d(x,y)$ on $\Gamma_H$ 
now measures the minimum number of single errors required to transform $x$ to $y$
\,or\, $y$ to $x$.

It is easily seen that every connected  simple graph $\Gamma$ can be
represented as an error graph where we can assume in addition  that
the  single error set consists of involutions (that is, partial maps $h$ defined on  suitable  subsets of $V$ such that $h^{-1}=h$). For if $c\!: E\rightarrow \{1...\chi\} \subseteq\mathbb{N} $ is an edge colouring of $\Gamma$ then each fiber $c^{-1}(i)$ with $i=1,...,\chi$ defines a natural involutionary error $h_{i}$ which is obtained by interchanging the two end vertices of any edge coloured by $i.$ In particular, every connected graph  $\Gamma$ is an error graph  with at most $\chi=\chi_{E}(\gG)$ errors where $\chi_{E}(\gG)$ is the edge-chromatic number of $\gG.$

 By Vizing's theorem~\cite{viz}\, this minimum number (over all $\Gamma$)
is equal to $k(\Gamma)+1$ where $k(\Gamma)$ is the maximum degree of
$\Gamma$, as in~(\ref{e008}). It is  a natural question to ask whether any
connected simple graph $\Gamma$ can be represented as an error graph $\gG_{H}$ for some error set $H$ of cardinality $k(\Gamma).$ The answer is  affirmative, see \cite{eros}, where it is also shown that the property (i) can in general not be replaced by a stronger
property $H=H^{-1}$ (meaning that $h^{-1}\in H$ if $h\in H$). 

In the examples discussed before,  the Hamming graph is an error graph when $V=F_q^n$  and  when  $H$ consists of the  $n(q-1)$ 
non-zero vectors $h\in F_q^n$ of Hamming weight 1,  with action given by  
$h(x)=h+x$ for  $x\in V$. Also the Johnson graph $J_w^n$ is of this
form when we view $V$ as the set of  all $w$-element subsets of $\{1,..,n\}$ and when $H$ is the set of all ${n\choose 2}$ transpositions  
$(i,j)$ interchanging $i$ and $j$ in $\{1,..,n\},$ in their natural permutational action 
on $V$ obtained by permuting  the coordinates of  vectors. In order to make sure that the single error property $h(x)\neq x$ holds for all vertices $x\in V_{h}$ one has to restrict the domain of $(i,j)$ to those sets which contain exactly one of $i$ and $j$.

Similarly, the insertion and deletion errors for finite sequences over an  alphabet $A$ can be 
described in this fashion as an infinite error graph. As vertex set we consider the set $V=A^0\cup A^1\cup A^2\cup...\cup A^n\,\cup\,...$\,\, of all finite words over $A.$  As single error set we take $H:=\{d_{1},\,d_{2},\,..,\,d_{m},\,...\}\cup \{i_{1}(a),\,i_{2}(a),\,...,\, i_{m}(a),\,...\quad :\, a\in A\}$ where $d_{m}$ deletes the $m^{th}$ entry in any word of length $\geq m$ while $i_{m}(a)$ inserts $a$ as the $m^{th}$ entry in any sequence of length $\geq m-1$. As expected, the usual graph metric is indeed the Levenshtein error distance~\cite{lev65} \,for sequences. Situations where the model of undirected single error graphs  is not applicable include 
asymmetric errors,  some further comments can be found in~\cite{eros}.

%However, for the sets $H$ of $n$ asymmetric errors on all
%$x=x_1,\ldots,x_n \in F_q^n$ each of which substitutes the symbol
%$x_i$ in x by the symbol $x_i+1$ or $x_i-1$, it is not true that
%$H=H^{-1}$ and this restriction in the definition of $\Gamma_H$ is
%weakened by so--called parallelogram property in \cite{eros}.

%{\scriptsize J: (1) We mentions these two graphs later on as Cayley
%graphs: It would be nice to mention also one general class of graphs
%of the shape $\gG_{H}$ which are not of the special form of Cayley
%graphs. (2) Also, we should mention error problems which can not be
%presented by \SE s. For instance, deletion and insertion errors are
%not of that kind, they can not be inverted uniquely. }

%{\scriptsize E: I added something concerning the previous J-remarks
%in the end of this section but I am sure that we have to write it
%more accurately.}

\section{Some Bounds for Regular Graphs}

For the remainder we assume that $\gG$ is a connected and regular graph on $v\geq 4$ vertices, with degree $2\leq k=k(\Gamma)$ and parameters
$\lambda=\lambda(\Gamma)$, $\mu=\mu(\Gamma)$. We have $0\le \lambda
\le k-1,$\, $1\le \mu \le k $ and the diameter of $\Gamma$ is $d(\Gamma)\ge 1.$

For some classes of regular and strongly regular graphs on $v$ vertices
we have $N(\Gamma,1)=o(v)$ as  $v\to \infty.$ The following  strongly regular graphs are well known, see Chapter~21~in~\cite{Wilson} or \cite{bcn}. The triangle graph $T(m)$ is strongly regular with
parameters $v=m(m-1)/2,$ $k=2(m-2),$ $\lambda=m-2,$ $\mu=4$ and
hence $N(T(m),1)=m.$ The lattice graph
$L_2(m)$ is strongly regular with  parameters $v=m^2,$ $k=2(m-1),$ $\lambda=m-2,$
$\mu=2$ and hence $N(L_2(m),1)=m.$ Meanwhile the 
Paley graphs $P(q)$ ($q$ a prime congruent to $1$ mod $4$) is strongly regular with parameters $v=q,$ $k=(q-1)/2,$ $\lambda=(q-5)/4,$
$\mu=(q-1)/4$ and hence $N(P(q),1)=(q+3)/4.$ The complement
of a strongly regular graph $\gG$ is also strongly regular (although not necessarily
connected). This complementary   graph $\overline {\Gamma}$ has parameters $v(\overline{\Gamma})= v,$ $k(\overline{\Gamma})=
v-k-1,$ $\lambda(\overline{\Gamma})=
v-2k-2+\mu,$ $\mu(\overline{\Gamma})= v-2k+\lambda,$
and hence $N(\overline{\Gamma},1)=v-2k+\max(\mu,\lambda).$

%{\scriptsize E: The graph $K'_v$ which is defined below isn't used
%in the text later. So it is not necessary to denote it by a special
%notation. Moreover, one can use $K^t_m$ instead of $K^{(t)}_m.$ But
%it is not so essential. One can rewrite the below part of the text
%as follows.}

Let $O^{t}_m=O_{m}\,\star O_{m}\,\star\,...\,\star\,O_{m}$ be the product of $t$ copies of the empty graph on $m$ vertices. This is  the complete $t$-partite graph with $v=tm,$ each part consisting of $m$ vertices and edges connecting vertices from different parts in all possible ways. If $t\geq 2$ this graph is connected and  strongly regular. 

The complete graph on $v$ vertices is denoted by $K_v.$ We recall that a
1-factor of a graph is a collection of disjoint edges covering all vertices (a complete matching of the vertices of $\Gamma$). When $v$ is even consider the graph obtained from $K_v$ by removing the edges of a
1-factor. This graph is strongly regular  with parameters $k=\mu=v-2$, $\lambda=v-4$ and coincides with $O^{t}_2$ with $t=\frac v2.$ Conversely, if $\Gamma$ is a regular of degree $k=v-2$ then $v$ is even and $\Gamma=O^{t}_2$ with $t=\frac v2.$ When $t=\frac v2$ then $N(O^{t}_2,1)=\lambda+2=v-2=\frac12(v+\lambda).$ More generally we have:

\begin{theorem}
\label{th001} \quad Let $\Gamma$ be a regular graph with $k\leq v-2.$  Then we have 
\begin{equation}
\label{e009} N(\Gamma,1)\le \frac 12 (v+\lambda)
\end{equation}
with equality if and only if $k-\lambda=v-k$ divides $v$ and 
$\Gamma$ is the strongly regular graph $O^{t}_{m}$ with  $m=k-\lambda$ and $v=tm.$  %$t=\frac {2k-\lambda}{k-\lambda}$ is an integer and $\Gamma$ is the strongly regular graph $O^{t}_{k-\lambda}.$
\end{theorem}

\pf By  \,(\ref{e007})\, we have $N(\Gamma,1)=\max\{\lambda+2,\,\mu\}.$ If $\Gamma=O^{t}_{m}$ with $v=tm$ then $k=v-m,$ $\mu=v-m$ and $\lambda=v-2m$\, so that $N(\Gamma,1)=\mu=\frac12(v+\lambda).$ For the converse assume first that $\lambda=k-1.$ In this case   \,(\ref{e007})\, implies that $N(\Gamma,1)=k+1$ which is not possible as $k+1$ is the cardinality of any single ball. Therefore $\lambda\leq k-2$ and from the assumptions in the theorem it follows  that $\lambda\leq v-4$ or $\frac12 \lambda + 2\leq \frac12 v.$   Hence $\lambda+2\leq \frac12(v+\lambda)$ 
with equality if and only if $\lambda=v-4.$ In the latter  case only $k=v-2$ is possible and so we have the situation already discussed, $\Gamma$ is $O^{t}_{m}$ with $m=k-\lambda=v-k=2$ and $v=2t.$

It is left to show that $\mu \le \frac 12 (v+\lambda)$ and to find the
conditions for  equality. For a $k$-regular graph $(V,\,E)$ and a vertex $x$ in $V$ we count the  number of edges between $S_1(x)$ and $S_2(x).$ This gives  
$$\sum_{y\in S_1(x)}(k-1-a_1(x,y))=\sum_{z\in S_2(x)}c_2(x,z),$$
see again the definitions in ({\ref{e004}). This gives $ k(k-1-\lambda)\le \mu k_2(x)$  and since $k_2(x)\le v-k-1$ we obtain  
\begin{equation} \label{e010}
k(k-1-\lambda)\le \mu k_2(x)\le \mu (v-k-1).
\end{equation}
 Since $1\le \mu\le k$ we get $k-1-\lambda\le v-k-1$ and hence $\mu\le k\le \frac 12(v+\lambda)$ as required.

If $\mu= k= \frac 12(v+\lambda)$ then we have  equalities in  (\ref{e010}). For a regular graph it is well-known that the inequalities in  (\ref{e010}) turn into equalities  if and only if the graph is strongly regular, see for instance {\sc Problem~21A}~in~\cite{Wilson}. So let $\Gamma$  be  strongly regular with $\mu=k.$ Then any pair of distinct and non-adjacent vertices have the same $k$ neighbours. It follows that $x=x'$ or $x$ is not adjacent to $x'$ defines an equivalence relation on the vertices of $\gG,$ with all equivalence  classes of size $m:=v-k.$ Hence $m$ divides $v=tm$ and $\Gamma=O^{t}_{m}.$
\hfill$\Box$

%Note that it follows from the proof that we have for any $k$-regular graph $\Gamma$ the inequality  \begin{equation} \label{e010} k(k-1-\lambda)\le \mu k_2(\Gamma). \end{equation}

\begin{theorem}
\label{th002} {\rm (Linear Programming Bound)} \quad Let $\Gamma$ be a regular graph of valency $k\geq 2.$ Then 
\begin{equation}
\label{e011} N_2(\Gamma,2)\ge\mu
\left(k-1-\frac 12(\mu
-1)(N(\Gamma,1)-2)\right)+2.
\end{equation}
\end{theorem}

We note that this rather general bound is  quite sharp, see the comment following Theorem~\ref{th005}. 

\medskip
\pf There are two vertices $x,x'\in V$ with $d(x,x')=2$ so that the set  $Y=\{y_1,..,y_{\mu}\}$ of all vertices at distance $1$ from both $x$ and $x'$ has $\mu\geq 1$ elements. %For $\mu=0$ the inequality \,(\ref{e011})\, holds.   
If $\mu=1$ then $y_{1}$ has $k-2$ neighbours other than $x$ and $x'.$ It follows that $ N_2(\Gamma,2)\ge |B_{2}(x)\cap B_{2}(x')|\geq 3+ k-2$ and so   \,(\ref{e011})\, holds. Hence we assume that $\mu\ge 2.$

\medskip
Let $U=\bigcup_{i=1}^{\mu}B_1(y_i)\setminus\{x,x'\}.$ We show that the number of elements in $U$ is  at least $\mu\left(k-1-\frac 12(\mu -1)(N(\Gamma,1)-2)\right).$  For $h=1,...,\mu$ let $U(h)$ be the vertices of $U$ which
belong to exactly $h$ of the sets $B_1(y_i),$ as $i=1,...,\mu$. 
In particular, $U=U(1)\cup...\cup U(\mu)$ is a partition and so 
$$|U|=\sum_{h=1}^{\mu}|U(h)|\,\,.$$
Next observe that the set $\left\{\left(z,\,B_{1}(y)\right)\,:\, y\in Y\mbox{\,\,and\,\,}z\in B_{1}(y)\cap U\,\right\}$ has cardinality $$\sum_{h=1}^{\mu}\,\,h|U(h)|=\mu(k-1)\,\,$$
 and the set $\left\{\left(z,\,\{B_{1}(y),\,B_{1}(y')\}\right)\,:\, y\neq y'\in Y \mbox{\,\,and\,\,}z\in B_{1}(y)\cap B_{1}(y')\cap U\,\right\}$ has cardinality 
$$\sum_{h=2}^{\mu}{h\choose 2}|U(h)|= \sum_{\{y,y'\}\subseteq Y,\,y\neq y'}\,\,(|B_1(y)\cap B_1(y')|-2)\leq {\mu\choose 2}\left(N(\Gamma,1)-2\right)\,.\spa$$
The last inequality holds as  $y\neq y'$ implies 
$|B_1(y)\cap B_1(y')|\le N(\Gamma,1).$ 

\medskip
Set $u_{h}:=|U(h)|$ for $h=1,...,\mu$ and $u:=|U|.$ To find a lower bound for $u$  we  minimize 
$$u-\mu(k-1)=-1u_{2}-2u_{3}-...-(\mu-1)u_{\mu}$$
for the non-negative integers $u_{2},\,...,u_{\mu}$ subject to the constraints
$$\sum_{h=2}^{\mu}\,\,hu_{h}\leq\mu(k-1)\,\,$$
and 
$$\sum_{h=2}^{\mu}{h\choose 2}u_{h}\leq {\mu\choose 2}\left(N(\Gamma,1)-2\right)\,.\spa$$

\medskip
Let $u^{*}\leq u-\mu(k-1)$ be  the required minimum. Then by the duality of linear programming, see for instance Section 7.5 in~\cite{Noble}, the value of   $u^{*}$ maximizes
$$-\mu(k-1)n_{1} \,\,-\,\,{\mu\choose 2}\left(N(\Gamma,1)-2\right)n_{2}$$ 
subject to $n_{1},\,n_{2}\geq 0$ and the dual constraints 
$$hn_{1}+{h\choose 2}n_{2}\geq h-1 \mbox{\,\,\, for\,\,\, }h=2,...,\mu\,\,.$$
Note that $n_{1}=0$ and $n_{2}=1$ satisfies the dual constraints  for all $\mu\geq 2$ and hence
$$u-\mu(k-1)\geq u^{*}\geq - {\mu\choose 2}\left(N(\Gamma,1)-2\right)\,\,.$$ Therefore 
$u\geq \mu\left(k-1-\frac{1}{2}(\mu-1)(N(\Gamma,1)-2)\right)$
as required.\dne

\medskip
Note for instance that  $N_2(\Gamma,2)\ge k+1$ when $\mu=1,$
$N_2(\Gamma,2)\ge 2k$ when $\mu=2$ and $N(\Gamma,1)=2$, and
$N_2(\Gamma,2)\ge 3k-4$ when $\mu=3$ and $N(\Gamma,1)=3$.

\begin{corollary}
\label{cor1} \,\,  Suppose that $\Gamma$ is a regular graph of valency $k$ with no triangle nor pentagons. If   $\mu \ge 2$ and $k\ge 1+{\mu \choose 2}$ then  
$$N_2(\Gamma,2)\ge N_1(\Gamma,2).$$
\end{corollary}

\pf   We have $\lambda =0$ since $\Gamma$ has no triangles so that  $N(\Gamma,1)=\mu$ by~(\ref{e007}). Similarly,  $N_1(\Gamma,2)=2k$ as $\Gamma$ contains no pentagons.  Using~(\ref{e011})\, we get
\begin{eqnarray}
N_2(\Gamma,2)-2k &\ge&\mu\left(k-1-\frac 12(\mu-1)(N(\Gamma,1)-2)\right)+2-2k  \nonumber \\
&=&\mu\left(k-1-\frac 12(\mu-1)(\mu-2)\right)+2-2k  \nonumber \\
&=&  (\mu-2)\left(k-1-{\mu\choose 2}\right)\ge 0\nonumber
\end{eqnarray}
and this completes the proof.
\hfill$\Box$

\bigskip

\section{Single error sets as group generators}

An important class of graphs associated to single error sets  is
obtained when the vertex set of the  graph are the elements of a
finite group. So we let  $G$ be a finite group and  consider the elements of $G=V$ as the vertices of the error graph $\gG=\gG_{H}$ for some error set $H.$ The neutral element of $G$ is denoted by $e=e_{G}$ and $1=\{e_{G}\}$ is the identity subgroup of $G.$
We suppose that  the single error set {\it is determined
as a subset $H$ of $G$} so that the action of errors on vertices is given by the group product. That is,  if $h\in H$ and $x\in G$ then $h(x):=xh^{-1}$.   In this situation  $H$  is a single error set  if and only if  $H$ does not contain $e_{G}$ and 

\vspace{-0.5cm}
\begin{enumerate}[(i)]
\item $H$ satisfies 
$H=H^{-1}(=\{\,h^{-1}\,\,:\,\,h\in H\,\})$, and
\item $H$ generates $G$ as a group.
\end{enumerate}

\vspace{-0.6cm}
The first condition is clear since there is some $g$ in $H$ with $g(h(x))=(xh^{-1})g^{-1}\\=x$ for a vertex $x$ in $V$ if and only if  $g=h^{-1}$ belongs to $H.$ The second condition is a restatement of the connectedness of the error graph. Note that we have set $h(x):=xh^{-1},$ rather than $h(x):=xh.$ This is advisable  so that  the multiplication of errors as elements of $G$ agrees with the co-cattenation of the corresponding maps, $(gh)(x)=x(gh)^{-1}=xh^{-1}g^{-1}=g\big(h(x)\big).$

As is well known, in this situation $\gG_{H}$ is the {\it undirected 
Cayley graph} on $G$ for the generating set $H,$ and $H$ 
is the {\it Cayley set} for $\gG_{H}$. Note conversely that every undirected Cayley graph can be viewed as a single error graph.

In the following we review some of the theory of  Cayley graphs from the 
viewpoint of single error graphs. Let $H$ 
be a Cayley set in the finite group $G$ with corresponding graph
$\gG_{H}=(V,\,E)$  on the vertex set $V=G$ and let ${\rm Aut}(\gG_{H})$ be the automorphism group of $\gG_{H}.$
 We consider two basic  kinds of
automorphisms of $\gG_{H}$. For each $g$ in $G$ the left-multiplication on $V,$ with  $g\!:\, x\mapsto gx$ for $x\in V,$   induces  an automorphism of $\gG_{H}$ since
$g\!:\,\{x,\,xh^{-1}\}\mapsto \{gx,\,gxh^{-1}\}$ maps edges to edges. If we think of $\{x,\,xh^{-1}\}$ as being labelled by $\overline{h}=\{x^{-1}(xh^{-1}),\,(hx^{-1})x\}=\{h^{-1},\,h\},$ the quotients of its end vertices, then  $\{gx,\,gxh^{-1}\}$ has the same label as   $\{x,\,xh^{-1}\}.$ Therefore   left-multiplication by elements of $G$ are automorphisms that preserves all  edge labels.

This action is transitive on vertices and only the identity element fixes 
any vertex. This is therefore the {\it regular action} of $G$ on itself. 
%By a theorem of Sabidussi~\cite{Sa58} 
This property characterizes Cayley graphs: $\gG$ is the Cayley graph 
of some group if and only if $\gG$ admits a  group of automorphisms  that acts regularly on its vertices, see for instance Chapter~6~in~\cite{Beineke}.  Note  however that  the graph usually does not determine the group.

We now describe graph automorphisms that change edge labels. Let $C$ be a group of automorphisms of $G$ as a group. For the action of $\gb\in C$ we write $\gb\!:\,x\mapsto \gb(x)$ and so $\gb(xy)=\gb(x)\gb(y)$ as $\gb$ is an automorphism of the group structure. We will also require that $C$ preserves $H,$ in the sense that $\gb(h)\in H$ for all $h\in H$ and $\gb\in C.$ Then $C$ is a group of automorphisms of $\Gamma_{H}$ since $\gb\!:\,\{x,\,xh^{-1}\}\mapsto \{\gb(x),\,\gb(xh^{-1})\}=\{\gb(x),\,\gb(x)\gb(h)^{-1}\}$ maps edges to edges, as   $\gb(h)^{-1}=\gb(h^{-1})\in H.$ Note that  the label of  $\gb(\{x,\,xh^{-1}\})$ now is $\overline{\gb(h)}.$ 

%Every $\ga\in A$ fixes $e=e_{G}\in G$ and permutes the neighbours of $e$ according to $h\mapsto \ga(h).$ Similarly, each $\ga\in A$ preserves the sphere $S_{i}(e)$ for each $i\geq 0$ and so $S_{i}(e)$ is a union of orbits for the action of $A$ on $V.$

The semi-direct  product $G\!\cdot\! C$ is the (abstract) group of all pairs $(g,\,\gb)$ with multiplication $(g',\gb')(g,\gb)=(g'\gb'(g),\,\gb'\gb).$ It acts on the graph as automorphisms by setting 
$$(g,\gb)\!:\,x\mapsto g\gb(x)\mbox{\quad for $x\in V.$}$$ This gives an injective group homomorphism 
from $G\!\cdot\! C$\, to ${\rm Aut}\,\Gamma_{H}$ so that we can regard $G\!\cdot\! C$ as a subgroup of ${\rm Aut}\,\Gamma_{H}.$ We collect these facts:

\begin{prop} \label{lem0002} \, Let $\Gamma_{H}$ be the error graph on the group $V=G$ with error set $H.$ Then the  left-multiplication of vertices by elements of $G$ forms a group of automorphisms of $\Gamma_{H}$ which acts regularly on the vertex set $V.$ If $C$ is a group of automorphisms of $G$ (as a group) such that $\gb(H)\subseteq H$ for all $\gb$ in $C$
then the semi-direct product  $G\!\cdot\! C$ is contained in the  automorphism group of $\Gamma_{H}.$ % be the automorphisms of $\Gamma_{H}$ defined via conjugation by elements from the normalizer $N=N_{G}(H)$ of $H$ in $G$.  Then the direct product $G\times C$ is a subgroup of  ${\rm Aut}\,\gG_{H}$ when we regard $(g,\,n)$ as to automorphism $x\mapsto gnxn^{-1}.$ Here  $G$ acts regularly on the vertices of $\Gamma_{H}$ by multiplication on the left, and $C$ fixes the vertices in the center $Z$ of $G$.  In particular, for any $z\in Z$ and any $i\geq 0$ the sphere $S_{i}(z)$ is a union of conjugacy classes of $N.$
\end{prop}

%We now describe automorphisms that change edge labels. Let $N:=N_{G}(H)=\{n\in G\,:\, n^{-1}hn\in H \mbox{ \,\,for all\, }h\in H\}$ be the normalizer in $G$ of $H,$ as a set. If $n$ is in $N$ then also right-multiplication by the inverse, $n\!:\,x\mapsto xn^{-1}$ for $x\in V$,  is an automorphism since $n\!:\,\{x,\,xh\}\mapsto \{xn^{-1},\,xhn^{-1}\}=\{xn^{-1},xn^{-1}(nhn^{-1})\}$ maps edges to edges, as $nhn^{-1}$ belongs to $H$ for all $h\in H$. Note, the edge $\{xn^{-1},\,xhn^{-1}\}$ now is labelled by $\overline{h^{*}}$ where $h^{*}=nhn^{-1}.$ Composing left-\, with  right-multiplication we see that  the conjugation map $x\mapsto nxn^{-1}$ is an automorphism of $\gG_{H}$ for all $n\in N.$ Let $C$ be the group of all such conjugation automorphisms induced by elements from $N.$ As $n\in N$ fixes all vertices of $\Gamma_{H}$ if and only if $n$ belongs to the center  $Z=\{z\in G\,:\, zg=gz \mbox{\,\, for all $g$ in $G$}\}$ of $G$ we see that $C$ is isomorphic to the quotient of $N/N\cap Z.$ %Note that the left- and right-multiplication actions commute, as a consequence of the associativity of multiplication.   The orbit  of  $x\in V$ under $C$ is the conjugacy class $x^{N}:=\{ nxn^{-1}\,:\,n\in N\,\}$. We note

A common example of this situation occurs when we consider conjugation by group elements. Let $b\in G.$ Then  {\it conjugation by $b$} is the  automorphisms 
$$x\mapsto  bxb^{-1}=:x^{b} \mbox{\quad for $x\in G$}$$ and $$x^{G}:=\{x^{b}\,\,:\,\,b\in G\}$$ is the {\it conjugacy class} of $x.$ In this case the error set $H$ is invariant under conjugation if and only if $H$ is a union of conjugacy classes. For $C$ we then take the group $C:=G/Z$ where $Z=Z(G)=\{b\in G\,:\, x^{b}=x \mbox{\,\, for all $x$ in $G$}\}$ is the {\it center } of $G.$ These are the {\it inner automorphisms } of $G.$ So in this situation $G\!\cdot\! C$ is a group of  automorphism of  $\Gamma_{H}.$ In Chapter 5 we will analyse this example further when $G$ is the symmetric group $\Sn$ on the set $\{1..n\}$ and when $H$ is the set of all transpositions on  $\{1..n\}.$ There we shall see that 
the full automorphism group  of the error graph can be larger than $G\!\cdot\! C,$ even if $C$ is the group of all automorphisms of $G$ as a group.

Another interesting example occurs when $\gG$ is the Hamming graph. Here $G$ is the vector space  $F_{q}^{n}$ where $F_{q}$ is the field of $q$ elements and $H$ is the set of all vectors of the shape $(0,..,0,a,0,..,0)$ with $a\neq 0$. Then $G$ acts on itself as a group of translations, that is, maps of the kind $g\!:x\mapsto g+x$ for all $x\in F_{q}^{n}.$   For $C$ we can take the {\it monomial subgroup } $C=(F_{q}^{\times})^{n}\cdot {\rm Sym}_{n}\subseteq{\rm GL}(n,q)$ acting naturally as linear maps  on $V.$ More precisely,  $C$ is the group of all $n\times n$ matrices with exactly one element from the multiplicative group $F^{\times}_{q}$ in each row and column. So here $F^{n}_{q}\cdot C$ is a group of affine linear maps on $F^{n}_{q}$ that acts naturally as automorphisms on the Hamming graph $\gG.$

\medskip
Considering again the general case we let $\Gamma_{H}=(G,\,E)$ be an error graph with error set $H.$ We have seen that any group automorphism $\gb$ fixing $H$ as a set induces an automorphism of $\Gamma_{H}.$ Evidently $\gb$ also fixes the identity element $e=e_{G}$ in $G.$ Assume therefore more generally that $C$ is a group of automorphisms of $\Gamma_{H}$ which fixes $e.$  For any $x\in V$$$x^{C}:=\{\gb(x)\,\,:\,\,\gb\in C\}$$ is the {\it orbit} of $x$ under $C.$ In order to analyze the parameters $k_{i}(x),$ $a_{i}(x,y),$ $b_{i}(x,y)$ and $c_{i}(x,y)$ note that $\Gamma_{H}$ is vertex transitive and therefore it suffices 
to consider the spheres with center $e_{G}$. Hence we abbreviate all parameters, writing $S_{i}$, $B_{i}$, $k_{i}=|S_{i}|$, $a_{i}(y)$, $b_{i}(y)$ and $c_{i}(y)$, suppressing the reference to $x=e_{G}$ in each case. In general these parameters still depend on $y$ although  automorphisms provide at least  for some form of regularity:

\begin{prop} \label{lem0003} \,
Let $\Gamma_{H}$ be the error graph on the group $V=G$ with error set $H$ and suppose that $C$ is a group of automorphisms of $\Gamma_{H}$ which fixes $e=e_{G}$. Then for each $i\geq 0$ the sphere $S_{i}=S_{i}(e)$ is a union of $C$-orbits. 

\vspace{-0.3cm}
Further, suppose that $y$ and $y'$ belong to the same $C$-orbit and that $r,\,t\geq 0.$ Then $|S_{r}\cap S_{t}(y)|=|S_{r}\cap S_{t}(y')|$ and $|B_{r}\cap B_{t}(y)|=|B_{r}\cap B_{t}(y')|.$ In particular, 
$a_{i}(y)=a_{i}(y')$,  $b_{i}(y)=b_{i}(y')$ and $c_{i}(y)=c_{i}(y')$ for all $i\geq 0.$ 
\end{prop}

\pf Let $y\in S_{i}$ and let $e,\,y_{1},\,...,\,y_{i}=y$ be a shortest path from $e$ to $y.$ If $\gb\in C$ then it is clear that $\gb(e)=e,\,\gb(y_{1}),\,...,\,\gb(y_{i})=\gb(y)$ is a shortest path from $e$ to $\gb(y).$ It follows that  $\gb(S_{i})=S_{i}$ is a union of $C$-orbits. Now suppose that $y'=\gb(y).$ Then
$\gb\left(S_{i}(y)\right)=S_{i}\left(\gb(y)\right)= S_{i}(y') $ and so $|S_{r}\cap S_{t}(y)|=
|\gb\left(S_{r}\cap S_{t}(y)\right)|=|\gb\left(S_{r}\right)\cap \gb\left(S_{t}(y)\right)|=|S_{r}\cap S_{t}(y')|.$ The remainder  follows immediately, including   the statement on  $a_{i},\,b_{i}$ and $c_{i}$ since these numbers are of the shape  $|S_{r}\cap S_{t}(y)|$ for particular choices of $r$ and $t.$ \dne

If $H$ is the single error set  of $\gG_{H}$ we set $H^{0}:=\{e_{G}\}$ and $H^{i}:=HH^{i-1}$ inductively for $i>0$.  Clearly, $\gG_{H}$ is regular of degree $k(\gG)=|H|$. If as before $S_{i}$ denotes the sphere of radius $i$ around $e=e_{G}$  then evidently $S_{1}=H^{1}=H,$ $S_{2}=H^{2}\setminus (H^{1}\cup H^{0})$ and more generally,  $$S_{i}=H^{i}\setminus (H^{i-1}\cup H^{i-2}\cup...\cup H^{1}\cup H^{0}).$$ The following is  easily shown and gives the value of  $N(\Gamma_H,1)$ by using~(\ref{e007}).

\begin{lemma}
\label{lem0005} \, In  an error graph  $\Gamma_H$ with error set $H$ we have 
$$
%\label{e018} 
\lambda(\Gamma_H) =\max_{x\in S_{1}}\mid \{(h,h')
\,:\,x=hh' \mbox{\,\,with\,\,}h,h'\in H \} \mid {\quad \rm  and}
$$$$
\label{e019} \mu(\Gamma_H) =\max_{x\in S_{2}}\mid \{(h,h')
\,:\,x=hh' \mbox{\,\,with\,\,}h,h'\in H \} \mid \,\,.
$$
\end{lemma}

\section{Permutations distorted by transpositional errors}

In the following we consider Cayley graphs when $G=\Sn$ is the symmetric group  acting on the set $\{1..n\}.$ Any subset $H$ of $G$ which generates $G$ with $e\not\in H$ and $H=H^{-1}$  is a Cayley set for $G.$ 

We express permutations in the usual cycle notation.  (Throughout the word `cycle' always refers to a particular kind of permutation, and never to a graph or subgraph.)  A {\it transposition} on  $\{1..n\}$ is a permutation of the shape $x=(i,j)$ with $1\leq i\neq j\leq n$ if we suppress the $1$-cycles of $x$. Particularly important graphs occur when $H=\{ (1,2),\, (2,3),...,\, (n-1,n)\}$ are the $n-1$ {\it Coxeter generators } of the symmetric group.  These form a minimal set of transpositions needed to generate $\Sn.$ This set corresponds to the fundamental reflections associated to a chamber for the $A$-type Dynkin diagram. The chambers give rise to a triangulation of the euclidean unit sphere in  $\mathbb{R}^{n-1}.$ In this situation the graph distance function 
$d(x,y)$ in $\gG_{H}$ is a discretized version of the geodetic distance on this sphere and presents the distance between two facets  in the triangulation of the sphere, see for instance the book~\cite{benson} \,of Grove and Benson on finite reflection groups. In this interpretation $B_{r}(x)$ is the `cap' of facets on the sphere at distance $\leq r$ from the facet $x$ and $N(\Gamma, r)$ is the number facets common to two such caps, with suitable distinct centers. 
Note also that here $d(e,\,-)$ evaluated for a single variable is the word length function in the corresponding Weyl group. This Cayley graph is of considerable importance in Lie theory and in many other parts of mathematics and  physics. For a recent treatment of its combinatorics we refer to \cite{Bjo05}.  We add that in computer science this graph is known as the {\it bubble--sort Cayley graph} and is used as a model for interconnection networks \cite{H97,LJD93}. 
 Various other Cayley graphs on $\Sn$ have been considered in the literature, we mention in particular Diaconis' book~\cite{Dia} where metrics on $\Sn$ more generally are discussed.

\medskip
By contrast we may consider the error graph on $\Sn$ when the single error set $H$ consists of {\it all } transpositions $(i,j)$ on $\{1..n\}.$ This clearly is a highly redundant system of generators, situated  at the other extreme to the case of the Coxeter elements in $\Sn$ which form a minimal generating set. 
In this situation a  single error $(i,\,j)$ transforms the vertex $x$ to its  neighbour $x(i,j)$ and all choices for  $1\leq i\neq j\leq n$ are admissible.   A graph $\gG_{H}$ of this type will be called a {\it transposition Cayley graph,} and these graphs are the subject of the remainder of the paper.

It may be useful to describe errors of this kind in a slightly more general setting. Let $A$ be a finite  alphabet with $|A|\geq 2$ and let $A^{n}$ be the set of all words of length $n$ over $A.$ Then the single transposition error $(i,\,j)$ on the coordinates of $A^{n}$ is the map $(i,\,j)\!:\, a=(a_{1},...,\,a_{i},\,...,\,a_{j},\,...,\,a_{n})\mapsto a^{(i,\,j)}= (a_{1},...,\,a_{j},\,...,\,a_{i},\,...,\,a_{n})$ with all other entries of $a$ unchanged. This gives rise to an error distance $d_{A}$ on $A^{n}$ where $d_{A}(a,\,b)$ is the least number of single transposition errors needed to transform $a$ to $b,$ if this is possible. In this case we must have $b=a^{g}$ for some $g\in \Sn$ and $d_{A}(a,\,b)\leq d(e_{G},\,g)$ where the latter denotes the distance in the transposition Cayley graph. (Observe that $d_{A}(a,\,a^{(i,\,j)})=0$ if and only if $a_{i}=a_{j}$ while $ d\left(e_{G},\,(i,\,j)\right)=1$ independently.) Note that this distance function  defines a graph on $A^{n}. $ Each component is  an error graph with involutory errors $(i,\,j)$ if we restrict the domain of the single error $(i,\,j)$ to the words $a$ in which $a_{i}\neq a_{j}.$ In this way the transposition error graph $\gG_{H}$ can be said to  control the transposition errors on $A^{n}.$   

In molecular biology transpositional errors are one  of the three known mechanisms in the mutation and evolution  of genetic information. The so-called {\it replication slippage} applied to a nucleotide sequence  is a process that results in some  strings of consecutive nucleotides being reversed or repeated in the sequence.
Such replication slippages usually recur and give rise to so-called microsatellites which  contain a high degree of information about the evolutionary process undergone by the nucleotide sequence in question, and often this happens in the non-coding part of the nucleotide sequence. For general information see Futuyma's book~\cite{Futu} on evolutionary biology as well as \cite{pev00} and \cite{sa02}.

Replication slippage is therefore a combinations of two kinds of errors on sequences, on the one hand the insertion-deletion process already mentioned at the end of Section 2 and the transpositional errors in the transposition Cayley graph on the other. It may be worth to mention that the other principal mutation mechanisms are {\it point mutations} referring to the replacement of one nucleotide by another, and {\it frame shifts } which are the insertion or deletion of a group of nucleotides. Both of these are therefore covered by the insertion-deletion process. 

Evidently any interval transposition or reversal (of a part of a nucleotide sequence) can be expressed as a product of single transpositional errors. However,  it should be interesting to introduce such products as new single errors, and to consider the resulting error graph on $\Sn.$ A second point  of interest should be to study the {\it resistance to transpositional errors:}\,  As the nucleotide alphabet consists of just four letters,  a single transpositional error is expressed only in a small proportion  of all possible words in $A^{n},$ leaving many others unchanged by that error.

\medskip
Returning to the general discussion of the transposition Cayley graph we note  the following conventions. Permutations in $\Sn$ are multiplied from the right so that $(xy)(j)=x(y(j))$ for all $x,\,y\in\Sn$ and $j\in \{1..n\}.$ If $x$ is written as a product of $h_{i}$ disjoint cycles of length $i$ for $1\leq i\leq n$
 then the {\it cycle type} of $x$ is denoted as ${\rm
ct}(x)=1^{h_{1}}2^{h_{2}}...\,n^{h_{n}}.$ Here it is essential to include $1$-cycles so that  $\sum_{i}\,\,i\,h_{i}=n.$ As is well-known, two permutations are conjugate to each other through  an element of $\Sn$ if and only if they have the same cycle type. Writing $G=\Sn$ therefore the conjugacy class $$(1^{h_{1}}2^{h_{2}}...\,n^{h_{n}})^{G}:=x^{G}=\{\,g^{-1}xg\,\,:\,\,g\in G\,\}$$ is the set of all permutations having the same cycle type as $x$.

We let $H=T:=\{\,(i,j)\in \Sn\,:\,1\leq i\neq j\leq n\}=(1^{n-2}\,2^{1})^{G}$ be the set of all
transpositions of $\{1..n\}$. Thus $\gG_{T}$  is the transposition Cayley
graph on $\Sn$ 
and will be denoted by $\Sn(T).$  The following collects some easily established facts.

\begin{lemma}
\label{lem0006} \, For $n \geq 3$ the transposition Cayley graph $\Sn(T)$
is a connected ${n\choose 2}$-regular graph of order $n!$
and diameter $n-1.$ It is $t$-partite for any $2\leq t\leq n.$
\end{lemma}

\pf The group $\Sn$ has order $n!$ and is generated by its ${n\choose 2}$ transpositions. Its diameter is at most $n-1$ since any permutation is a product 
of at most $n-1$ transpositions. 
On the other hand, an $n$-cycle can not be written in terms of fewer than $n-1$ transpositions. No two  elements in the same sphere $S_{i}$ could  be adjacent to each other  as they have the same determinant $(-1)^{i}.$ Hence $S_{0},$ $S_{1},$ $S_{2},$ ..., $S_{n-1}$ is a partition into $n$ parts from which a $t$-partition can be obtained  for any $t\leq n$.  \dne

For the product of a permutation with a transposition  the following simple rule is
essential.  If $x=(i_{1},..,\,i_{k})(j_{1},..,\,j_{\ell})$ consists of two disjoint cycles
and if $t=(i,j)$ interchanges elements from different cycles, say $i_{1}=i$ and
$j_{1}=j$ without loss of generality as the cycles are determined
only up to cyclic reordering, then
\begin{equation}
\label{e020}
xt=(i_{1},j_{2},j_{3},..,j_{\ell},j_{1},i_{2},i_{3},..,i_{k})=:s
\end{equation}  is a single cycle obtained by joining up the two cycles of $x$. Conversely, upon multiplying this equation again by $t$, we see that
multiplying the single cycle $s$ by a transposition of some two
elements from that cycle gives $x=xtt=st$, hence splitting that
single cycle into two cycles. Therefore multiplying any permutation $x$ by a transposition results in a permutation which either joins up two cycles of $x$ or splits one cycle of $x$ into two, with no other changes.

Following the earlier convention  whereby $S_{i}=S_{i}(e)$ we
have that  $H=T=S_{1}$ consists of all transpositions, $S_{2}$
consists of all $3$-cycles $(i,j,k)$ and all {\it double transpositions}
$(i,j)(k,\ell)$ with  $i,\,j,\,k,\,\ell$ distinct,  and so on. As multiplication by a transposition
increases or decreases the number of cycles by one\, it follows by
induction  that $S_{i}$ consists of all permutations expressible as
a product of $n-i$ disjoint cycles, counting also all $1$-cycles.

The path distance between  two permutations $x$ and $y$ is the least
number $d$ of transpositions $t_{i}$ such that $xt_{1}...t_{d}=y$.
Equivalently $d$ is the least number of transpositions needed to write
$x^{-1}y$ and also equal to the number of bisections and gluings
needed to transform the cycles of $x$ into those of $y$. The number
of distinct paths from $x$ to $y$ is equal to the number of paths from
$e$ to $x^{-1}y$ and about these the following theorem gives
complete information. It is based on Ore's theorem on the number of
trees with $n$ labeled vertices, see also Theorem 2 \,in~\cite{langlands}.

\begin{theorem}\cite{den}\quad 
\label{th003} Suppose that  $x$ has cycle type ${\rm
ct}(x)=1^{h_{1}}2^{h_{2}}...n^{h_{n}}$ and let $1\leq i \leq n-1$ be such that $\sum_{j=1}^{n}\,h_{j}=n-i$ is the number of cycles in $x.$ Then the number of distinct
ways to express $x$ as a product of $i$ transpositions is  equal to
\begin{equation}
\label{e021}
i!\prod_{j=1}^n\left(\frac{j^{j-2}}{(j-1)!}\right)^{h_j}.
\end{equation}
\end{theorem}

By the discussion above $x$ cannot be written in fewer than $i$ transpositions.   The special case $i=n-1$ and $h_n=1$ means that each of the $(n-1)!$
cycles of length $n$ has $n^{n-2}$ different representations as a
product of $n-1$ transpositions. This number coincides with the
number of trees with $n$ labelled vertices, see also Section 5.3 in \,Stanley
\cite{stanley2}.

Let $1\leq i \leq n-1.$ %and $1\leq j \leq n.$  
If $y\in S_{i}$  has cycle type ${\rm ct}(y)=1^{h_{1}}2^{h_{2}}..\,n^{h_{n}}$ and consists of $\sum_{j=1}^{n}\, h_{j}=n-i$ cycles then $y$ is a product of $i$ transpositions. As
$\det y=(-1)^{i}$ we must have $a_{i}(y)=0$. As a single cycle of
length $j$ can be split into two cycles as in~(\ref{e020}) in ${
j\choose 2}$ different ways, we have $c_{i}(y)=\sum_{j=1}^{n}\,
{j\choose 2}h_{j}$. From $\sum_{j=1}^{n}\, jh_{j}=n$ it follows
that $c_{i}(y)=\sum_{j=1}^n\,{j \choose 2}h_j=\frac 12
\left(\sum_{j=1}^n\,j^2h_j-n\right)$. If we regard $y$ as an element of ${\rm Sym}_{m}$ with $m>n$ then it is clear from~(\ref{e020}) that $c_{i}(y)$ is independent of $n$. Finally, $b_{i}(y)={n\choose
2}-c_{i}(y)$. We collect these facts:

\begin{lemma}
\label{lem0007} \, In $\Sn(T)$ the set $S_{i}$, where  $1\leq i
\leq n-1,$ consists of all  permutations 
of $\{1..n\}$ which are composed of  exactly $n-i$ disjoint cycles, including $1$-cycles.

\vspace{-0.3cm}
If $y\in S_i$ has cycle type  ${\rm
ct}(y)=1^{h_{1}}2^{h_{2}}..\,n^{h_{n}}$ then
$$ c_i(y)=\frac 12 \left(\sum_{j=1}^nj^2h_j-n\right),$$
$$ a_i(y)=0\mbox{\quad and}$$
$$ b_i(y)=\frac 12 \left(n^2-\sum_{j=1}^nj^2h_j\right).$$
If $y$ is regarded as an element in ${\rm Sym}_{m}$ with $m>n$ then only $b_{i}(y)$ depends on $n.$ 
\end{lemma}

{\sc Remark:}\quad Loosely speaking,  if  $y$ belongs to $S_{i}$ we can think of $c_{i}(y)$ as the `downward' degree of $y,$  namely the number of neighbours of $y$ in the next lower sphere $S_{i-1}.$ The fact that this degree is independent of $n$ will be used later on. Similarly ${n\choose 2}-c_{i}(y)$ is the `upward' degree of $y.$
The transposition Cayley graphs are not distance-regular and they  illustrate the fact that the up-\, and downward degrees are not constant for elements in the same sphere. This can be seen already in ${\rm Sym}_{4}(T).$ If $y$ in $S_{2}$ is a $3$-cycle then $c_{2}(y)=3$ according to the three choice of a transposition splitting the $3$-cycle. On the other hand,  if   $y=(1,2)(3,4)$ in $S_{2}$ is a double transposition then $c_{2}(y)=2$ as there are just two ways to split one of the two cycles. This is true for any  $n\geq 4.$ %On the other hand, any two elements $y$ and $y'$ of the same cycle type are in the same conjugacy, and so also Proposition~\ref{lem0003} implies that $a_{i}(y)=a_{i}(y')$ and $b_{i}(y)=b_{i}(y'),$ as we shall see now.

\medskip
Next we discuss the automorphism group of the transposition Cayley graph. As before let $G=\Sn=V$ and set $\Gamma=\Sn(T).$  Let $(a,b)$ be an element of the direct product $G\times G$. Then $(a,b)\!: x\mapsto axb^{-1}$ for $x \in V$ is an automorphism of $\gG$ since for any transposition $t$ we have $xt\mapsto axtb^{-1}=(axb^{-1})(btb^{-1})$ in which $btb^{-1}$ again is a transposition. Note that only the identity of $G\times G$ fixes all vertices  since $axb^{-1}=x$ for all $x\in V$ implies that $a=b $ and hence that $a\in Z(G)=1.$ This implies that we can view $G\times G$ as a subgroup of ${\rm Aut}(\gG).$ Recall the discussion in Section 4. If we let $C$ be the group of conjugation automorphisms,  $x\mapsto x^{b}$ for all  $b\in G,$ then  $C$ is the {\it  diagonal subgroup } $\{(b,b)\,:\,b\in G\}\subseteq G\times G.$ Furthermore, we have $G\times G=G\!\cdot\!C$ as subgroups of ${\rm Aut}(\gG).$

A further automorphism of $\Gamma$ comes from the  inversion map $$\imath\!: x \leftrightarrow  x^{-1}\mbox{\quad for $x\in V.$}$$ While $\imath$ is not an automorphism of the group it  is an automorphism of the graph. For if $\{x,y\}$ is an edge with  $y=xt$  and $t$ a transposition  then $y^{-1}=x^{-1}(yty^{-1})$ where  $yty^{-1}$ is a transposition and so  $\{y^{-1},x^{-1}\}$ is an edge. Since $\left(\imath(a,b)\imath\right)(x)=(ax^{-1}b^{-1})^{-1}=(b,a)(x)$ for all $x\in V$ we see that  $\imath$ normalizes $G\times G$ by interchanging the two direct factors. This shows that  the semi-direct product
$(\Sn\times \Sn)\cdot\langle \imath\rangle$ is contained in $ {\rm Aut}(\Sn(T))$.

\begin{theorem}
\label{con1} \quad For $n\geq 3$ the full automorphism group of $\Sn(T)$ is the semi-direct product  $(\Sn\times \Sn)\cdot C_{2}$ where $C_{2} =\langle \imath\rangle$ is the group of order $2$ obtained by inverting  the elements in $V=\Sn$. 
\end{theorem}

%The factors $G\times 1$ and $1\times G$ inside ${\rm Aut}(\Sn(T))$  are the left- and right-regular representations of $G$ which are interchanged by inversion in the group. From the argument above it is clear that both regular representations are contained in the automorphism group whenever the Cayley set is normal in the group. 

\pf  As before set $G=\Sn,$ $\gG=\Sn(T)$ and let $A$ be the group of all automorphisms of $\gG.$ When  $n=3$ when $\gG$ is the complete bipartite graph $K_{3,3}$ and in this case the statement can be checked directly from the description of the action of $({\rm Sym}_{3}\times {\rm Sym}_{3})\cdot\langle \imath\rangle$ on $\gG$. 

 Now suppose that $n>3$ and let $\ga''\in A.$ As $G$ acts vertex transitively by left-multiplication we select $g''\in (G\times 1)\subseteq (G\times G)$ such that $\ga':=g''\ga''$ fixes $e_{G}$. This implies that $\ga'$ fixes $S_{r}$ as a set, for all $r\geq 1$, see Proposition~\ref{lem0003}. 
%In the first part we show that an automorphism of $\Sn(T)$ fixing all points in $B_{1}(e_{G})$ and one $3$-cycle in $S_{2}$ must be the identity. In the second part we conclude the proof by showing that an arbitrary automorphism can be change by an element in $\Sn\cdot(\Sn\cdot C_{2})$ so that it fixes all the vertices just mentioned. 
Furthermore, $\ga'$ fixes each of the two conjugacy classes $(1^{n-3}3^{1})^{G}$ and  $(1^{n-4}2^{2})^{G}$ in $S_{2}$ since $c_{2}=3$ on the first class while $c_{2}=2$ on the second class, see  the remark following Lemma~\ref{lem0007}.

For $1\leq i\leq n$ let  the {\it pencil\,}  $P_{i}$ be the set $P_{i}=\{\,(i,j)\in S_{1}\,:\,1\leq j\leq n\mbox{\quad and \quad} i\neq j\}.$ Then the following holds: any pair $x\neq y\in P_{i}$ has exactly two joint neighbours in   $(1^{n-3}3^{1})^{G},$\, and $P_{i}$ is  a maximal subset of $S_{1}$ with this property. Conversely, any set of $n-1$ vertices in $S_{1}$ satisfying  this property is a pencil. Since $(1^{n-3}3^{1})^{G}$ is invariant under $\ga'$ we see that $\ga'(P_{i})$ again is a pencil. 
On the other hand, the diagonal element $(g,g)\in G\times G$ satisfies $(g,g)(i,j)=g\cdot(i,j)\cdot g^{-1}=(g(i),g(j))$  so that  $(g,g)(P_{i})$ is the pencil $P_{g(i)}$. This means that the diagonal group  induces the full symmetric group on pencils while fixing $e_{G}$. In particular,  we can find some $g'=(g,g)\in G\times G$ such that $\ga:=g'\ga'$  fixes each pencil as a set, in addition  to the vertex $e_{G}.$ Let $x=(i,j)$ be an element of $S_{1}$. Then $\{x\}=P_{i}\cap P_{j}$ so that  $\{\ga(x)\}=\ga(P_{i})\cap \ga(P_{j})=\{x\}.$ Hence $\ga\,$ fixes all elements in $B_{1}(e_{G})$ pointwise.  

Note that $(1,2),\,(1,3)$ and $(2,3)$ are pairwise joined to two elements in $S_{2}$, and no others, namely $x=(1,2,3)$ and $y=(1,3,2)$. Thus $\ga$ fixes $\{x,\,y\}$ as a set and if \,$\imath$\, denotes the inversion automorphism mentioned before, then either  
$\ga$ or $\imath\ga$ fixes all of $B_{1}(e_{G})\cup\{ (1,2,3)\}$ pointwise.  By the following lemma either $\,\imath g'g''\ga''$\, or  $g'g''\ga''$\, is the identity automorphism of $\gG$ and so  $\ga''=g''^{-1}g'^{-1}\imath\,$ or $\ga''=g''^{-1}g'^{-1}\,$ belongs to $(\Sn\times \Sn)\,\cdot\, C_{2}.$ \dne

\begin{lemma}
\label{lem0007a} \, For $n\geq 3$ only the identity automorphism of $\Sn(T)$ fixes 
every  vertex in $B_{1}(e_{G})\cup\{ (1,2,3)\}.$
\end{lemma}

\pf This is evident for $n=3.$ Suppose therefore that  $n\geq 4$ %and that  $G,$ $\gG$ and $A$ are as before. Suppose that 
and that $\ga$ is an automorphism fixing  
every  vertex in $B_{1}(e_{G})\cup\{ (1,2,3)\}.$

Then each  double transposition in $(1^{n-4} 2^{2})^{G}$ is fixed by $\ga$ as these elements have exactly two neighbours in $S_{1},$ with no two double transpositions having the same $S_{1}$-neighbours. The elements in  $(1^{n-3}, 3^{1})^{G}$ fall into pairs $[(i,j,k),\,(i,k,j)]$ of $3$--cycles, each pairwise linked to the three fixed elements $(i,j),$ $(j,k)$  and $(i,k)$ in $S_{1}.$ Therefore $\ga$ either fixes or interchanges the members in each pair. We show that $\ga$ fixes these elements and hence is the identity on $S_{2}$.

Evidently $(1,2,3)$ and $(1,3,2)$ are both fixed. Hence look at the three pairs $[(1,4,2), (1,2,4)]$,\, $[(1,3,4), (1,4,3)]$ and $[(2,3,4), (2,4,3)]. $ As can be calculated, the six $4$-cycles in $S_{3}$ involving $1,\,2,\,3$ and $4$ are partitioned into two sets $X$, all connected to $(1,2,3), $ and $Y,$ all connected to  $(1,3,2).$ The sets $X$ and $Y$ are therefore fixed by $\ga$ as sets. It turns out that $(1,4,2)$ is linked to two vertices in $X$ while $(1,2,4)$ is linked to two vertices in $Y$. This means that $(1,4,2)$ and  $(1,2,4)$ are each fixed by $\ga$. 
The same argument extends to all other $3$--cycles. Hence $B_{2}(e_{G})$ is fixed pointwise.  For the remainder the argument becomes more homogeneous. Suppose that $x$ and $y=\ga(x)$ are in $S_{r}$ with $ r>2.$ By induction we can assume that $\ga$ fixes all vertices in $S_{r-1}$ and this means that $x$ and $y$ have the same neighbours $N(x)=N(y)$  in  $S_{r-1}.$ We claim that this forces $x=y.$ The elements in $N(x)$ are obtained by 'splitting' any cycles appearing  in $x$ into two cycles in all possible ways, see (\ref{e020}). In particular, $x$ and $y$ have the same orbits on $\{1..n\}$  and if there are at least two orbits of length $>1$ then $N(x)=N(y)$ forces $x=y.$ In the remaining case $x$ and $y$ consist of a single cycle of length $\ell\geq 4$ with all other vertices fixed. It is easy to see that $\ell>3$ and  $N(x)=N(y)$ again forces $x=y.$ \dne }

\bigskip
\section{Distance statistics in the transposition graph}

Let  $S_{i}$ be the sphere of radius $i\leq n-1$ and centre  $e_{G}$ in the transposition Cayley graph $\Sn(T).$ Then $S_{i}$ is  a union of $\Sn$-conjugacy classes and  the parameters $a_{i}(y)$, $c_{i}(y)$ and $b_{i}(y)$ are constant on these classes, for all $1\leq i \leq n-1.$ It will  be useful to set $s_{i}(n):=|S_{i}|,$ and in more customary symbols, $c(n,n-i):=|S_{i}|.$

Then  $c(n,n-i)$ is the number of permutations in $\Sn$ having
$n-i$ cycles, for $1\leq i \leq n-1,$ and these  are the   {\it  signless Stirling numbers of the first kind, } see for instance Chapter~1.3\, in Stanley's book~\cite{stanley2}. 
We have  $c(n,n)=1$, $c(n,n-1)={n\choose 2}$, $c(n,n-2)=2{n
\choose 3}+3{n \choose 4},$ $c(n,n-3)=3{n\choose 4}+20{n\choose 5}+15{n\choose 6}$ and so on, up to $c(n,1)=(n-1)!$. The generating function of  $c(n,m)$ satisfies $$g(t):=\sum_{m=1}^nc(n,m)t^m=t(t+1)\cdots (t+n-1)\,\,$$ and from this we get the product form $$\Pi_{\Sn(T)}=t^{n}g(t^{-1})=(1+t)(1+2t)\cdots\left(1+(n-1)t\right)$$ for the Poincar\'{e} polynomial \,(\ref{poin})\, of $\Sn(T).$ From the definition it is clear that $s_{i}(n)$ is a polynomial in $n$ when $i$ is  fixed. The leading term counts the number of permutations of cycle type $1^{n-2i}\,2^{i}$ and so we note:

\begin{lemma}\label{srn} \, If $i$ is fixed and $n\geq 2i$ then $s_{i}(n)$ is a polynomial in $n$ of degree $2i. $ Its leading term is the leading term of $\frac{1}{i!}{n\choose 2} {n-2\choose 2}\cdots {n-2i+2\choose 2}$ and is equal to $\frac{1}{i!2^{i}}n^{2i}.$
\end{lemma}

\medskip
For $y\in S_{i}$ with cycle type ${\rm
ct}(y)=1^{h_{1}}2^{h_{2}}...\,n^{h_{n}}$ let
as before $(1^{h_{1}}2^{h_{2}}...\,n^{h_{n}})^{G}=y^{G}$ be the conjugacy
class of $y$. Then
$$|(1^{h_{1}}2^{h_{2}}...\,n^{h_{n}})^{G}|= \frac {n!}{1^{h_1}h_1!2^{h_2}h_2!\cdots
n^{h_n}h_n!}\,\,\,,$$ and 

\begin{equation}
\label{e022} S_i=\bigcup_{h_1+h_2+\cdots +h_n=n-i}
(1^{h_{1}}\,2^{h_{2}}\,...\,n^{h_{n}})^{G}\,\,,
\end{equation}

see again Chapter~1.3\, in~\cite{stanley2}. Omitting  cycle types of multiplicity $0$ 
we therefore  have $S_{1}=
(1^{n-2}\,2^{1})^{G},$\, $S_{2}=  (1^{n-3}\,3^{1})^{G}\,\cup
(1^{n-4}\,2^{2})^{G}$\, and so on. For small values of $r$ one can compute $N(\Sn(T),r)$ easily from this information. 

\medskip
\subsection{The value of $N(\Sn(T),r)$ for $r\leq 3$}

As we have observed, $\Sn(T)$ is not distance-regular and as a consequence it is not straightforward to determine the value $N(\Sn(T),r)$ for general $r.$ We begin to evaluate $N(\Sn(T),r)$ for $r\leq 3$ when closed formulae can be  obtained. 

\begin{theorem}
\label{lem0008} \quad For $n\ge 3$ we have 
\begin{equation}
\label{e023} N(\Sn(T),1)=3\,.
\end{equation}
\end{theorem}

\pf From Lemma~\ref{lem0007} we have $\lambda(\Sn(T))=0$ since
$a_1(z)=0$ for  $z \in S_1$ and moreover $c_2(y)=3$ if $y$ has cycle type $ct(y)=1^{n-3}\,3^1$ and
$c_2(y)=2$ if $y$ has cycle type $ct(y)=1^{n-4}\,2^2.$ Therefore,
from~(\ref{e004}) we have $\mu(\Sn(T))=3$ and by~(\ref{e007}) we
get~(\ref{e023}). \dne

\begin{theorem}
\label{th005} \quad For $n\ge 5$ we have 
\begin{equation}
\label{e024} N(\Sn(T),2)=N_2(\Sn(T),2)=\frac 32(n+1)(n-2).
\end{equation}
\end{theorem}

{\sc Remark:\,\,}  From this result one can see that the bound in Theorem~\ref{th002} is indeed very good: Working out the parameters for the transposition Cayley graph gives the bound $N(\Sn(T),2)\geq N_2(\Sn(T),2)\geq \frac 32(n+1)(n-2)-1$
from Theorem~\ref{th002}.

\pf By vertex transitivity it suffices to compute $|B_{2}\cap B_{2}(y)|$ with$B_{2}=B_{2}(e).$ This quantity depends only on the conjugacy class to which $y$ belongs, this is a consequence of  Proposition~\ref{lem0003} \,and Theorem~\ref{con1}.
Therefore we need  to consider the number $N(y)$  of all 
vertices in $B_2$ which are at distance $\leq 2$ from a given vertex $y \in S_i$ when $i$ runs from $1$ to $4.$  By~(\ref{e022})
 we have 
 
 \vspace{-0.6cm}
 $$S_4=(1^{n-5}\,5^1)^G\cup\,(1^{n-6}\,2^1\,4^1)^G\cup\,(1^{n-6}\,3^2)^G\cup\,
(1^{n-7}\,2^2\,3^1)^G\cup\,(1^{n-8}\,2^4)^G, $$
$$S_3=(1^{n-4}\,4^1)^G\cup\,(1^{n-5}\,2^1\,3^1)^G\cup\,(1^{n-6}\,2^3)^G\,$$

and so on. The numbers $N(y)$ are presented in Table 1. The row index is the conjugacy class which contains $y$ while the column  index is  the conjugacy classes contained in $B_{2}.$ The value of $N(y)$ is worked out using (\ref{e020}).

\begin{center}
\begin{tabular}{|l||c|c||c||c|}
\hline
\mbox{\,\,\,\,\,\,\,\,\,\,\,N(y) }& $(1^{n-3}\,3^1)^G$ & $(1^{n-4}\,2^2)^G$ & $(1^{n-2}\,2^1)^G$ & $(1^n)^G$ \\
\hline \hline
$(1^{n-5}\,5^1)^G$ & 10 & 10 & 0 & 0\\
\hline
$(1^{n-6}\,2^1\,4^1)^G$ & 4 & 6 & 0 & 0\\
\hline
$(1^{n-6}\,3^2)^G$ & 2 & 9 & 0 & 0\\
\hline
$(1^{n-7}\,2^2\,3^1)^G$ & 1 & 7 & 0 & 0\\
\hline
$(1^{n-8}\,2^4)^G$ & 0 & 6 & 0 & 0\\
\hline \hline
$(1^{n-4}\,4^1)^G$ & 4 & 2 & 6 & 0\\
\hline
$(1^{n-5}\,2^1\,3^1)^G$ & 1 & 3 & 4 & 0\\
\hline
$(1^{n-6}\,2^3)^G$ & 0 & 3 & 3 & 0\\
\hline \hline
$(1^{n-3}\,3^1)^G$ & $6(n-3)+2$ & $3{n-2\choose 2}$ & 3 & 1\\
\hline
$(1^{n-4}\,2^2)^G$ & $4(n-2)$ & $2{n-2\choose 2}-1$ & 2 & 1\\
\hline \hline $(1^{n-2}\,2^1)^G$ & $2(n-2)$ & ${n-2\choose 2}$ &
${n\choose 2}$ & 1\\
\hline
\hline

\end{tabular}\\
\vspace{0.5cm}
{\sc Table 1}
\end{center}

When we consider the corresponding rows in the table we get $N_4(\Sn(T),2)=20$ when
$n\ge 5,$ $N_3(\Sn(T),2)=12$ when $n\ge 4,$\, $N_2(\Sn(T),2)=\frac 32
(n+1)(n-2)$ and $N_1(\Sn(T),2)=(n-1)n$ for all $n\ge 3$. This  proves
the theorem due to~(\ref{e005}). \dne

\medskip
To estimate the number of vertices in  $|B_{r}\cap B_{r}(y)|$ for an arbitrary $y$ we consider the paths  $t_{1}t_{2}\cdots t_{r^{*}}$ with $r^{*}\leq r$ starting at $y$ and leading to a vertex $z=yt_{1}t_{2}\cdots t_{r^{*}}$ belonging to  $B_{r}.$ We say that this path has a  \,{\it descent} \,  at step $k<r^{*}$ if $yt_{1}t_{2}\cdots t_{k-1}\in S_{s}$ for some $s$ while $yt_{1}t_{2}\cdots t_{k-1}t_{k}\in S_{s-1}.$ The number of ways to continue the path at  $yt_{1}t_{2}\cdots t_{k-1}$ by a descent is the downward degree $c(e_{G},\,yt_{1}t_{2}\cdots t_{k-1})$ of \,(\ref{e00001})\, which by Lemma~\ref{lem0007}\, is independent of $n$.
Similarly, we say  that the path has an \,{\it ascent}\, at step $k$ if $yt_{1}t_{2}\cdots t_{k-1}\in S_{s}$  while $yt_{1}t_{2}\cdots t_{k-1}t_{k}\in S_{s+1}.$ In this case 
the number of choices to continue the path at  $yt_{1}t_{2}\cdots t_{k-1}$ by an ascent is the upward degree $b(e_{G},\,yt_{1}t_{2}\cdots t_{k-1})$ \, which by Lemma~\ref{lem0007}\, is of order $n^{2}.$ Hence

\begin{lemma}
\label{lem00006} \, The number of vertices $z=yt_{1}t_{2}\cdots t_{r^{*}}$ reachable from $y$ on a path with $a$ ascents is at most $k_{y}n^{2a}$ where $k_{y}$ is some constant independent of $n.$
\end{lemma}

Let $E_{i,i+1}$ be the set of edges joining a vertex in $S_{i}$ to one in $S_{i+1}.$ As $\Sn(T)$ is $k$-regular with $k=|S_{1}|$ and as $a_{i}(z)=0$ for all $z\in S_{i}$, see Lemma~\ref{lem0007},  we have $|E_{i-1,i}|+|E_{i,i+1}|=k\cdot|S_{i}|$  and hence
\begin{equation}
\label{eij1}
|E_{r-1,r}|=k\cdot (|S_{r-1}|-|S_{r-2}|+|S_{r-3}|-...+(-1)^{r-1}|S_{0}|)\,\,
\end{equation}
for all $r.$ For the transposition $y\in T$ let $E_{r-1,r}(y)$ be the set of all edges in $E_{r-1,r}$ of the form $\{z,zy\}$. (These are the edges in $E_{r-1,r}$ that  are labelled by $y$.) Evidently all automorphisms of $\Sn(T)$ fixing $e_{G}$ permute $E_{r-1,r}$ as a set and since the conjugation action is transitive on $T$ every edge label must appear an equal number of times in each orbit. Hence
\begin{equation} 
\label{eij2}
|E_{r-1,r}(y)|=|S_{r-1}|-|S_{r-2}|+|S_{r-3}|-...+(-1)^{r-1}|S_{0}|\,\,
\end{equation}
for all $r.$ If $y=(j_{1},j_{2})$ then the end vertices $v^{-}\!\in S_{r-1}$ and $v^{+}\!\in S_{r}$ of $\{v^{-},v^{+}\}\in E_{r-1,r}(y)$ are composed of cycles in which $j_{1}$ and $j_{2}$ occur in different, respectively the same, cycle(s). Hence (\ref{eij2}) gives the number of permutations in $S_{r-1}$ with $j_{1},\,j_{2}$ in different cycles and, at the same time,  the number of permutations in $S_{r}$ with $j_{1},\,j_{2}$ in the same cycle.

\begin{theorem}
\label{th006} \quad Let $\Gamma=\Sn(T)$ be the transposition Cayley graph and suppose that  $n\geq 4.$ Then  we have 
%\vspace{-0.6cm}
\begin{eqnarray}
\label{e0024} 
(i)\,\,\,\quad N_{1}(\gG,\,3)&=&2|S_{0}|+2|S_{2}| \,\,\,\,\,{\rm and} \nonumber\\
(ii)\,\,\quad N_{2}(\gG,\,3)&=&|S_{0}|+|S_{1}|+|S_{2}|+(n+2)(n-3)+\nonumber\\&+&24{n-3\choose 2}+22{n-3\choose 3}+6{n-3\choose 4}\,.\nonumber
\end{eqnarray}
Furthermore we have $N(\gG,\,3)=N_{2}(\gG,\,3)$ for all  $n\geq 16.$ 
\end{theorem}

\pf We need to compute $|B_{3}\cap B_{3}(y)|$ when $e=e_{G}$ and $y$ have distance $d(e,y)\leq 6$ from each other. When $d(e,y)=5$ or $6$ then a  path of length $\leq 3$ from $y$ to a vertex in $B_{3}$ can not have any  ascents. Therefore the number of such vertices is independent of $n$ by Lemma~\ref{lem00006}. (The same phenomenon can be seen in the upper part of Table 1.) 
When $d(e,y)=3$ or $4$ then the corresponding paths have at most one ascent so that 
$|B_{3}\cap B_{3}(y)|$ is of order at most $n^{2}.$ When $d(e,y)=1$ or $2$ then $|B_{3}\cap B_{3}(y)|$ is of order at least $n^{4}$ as we will show now.  It follows that the cases $3\leq d(e,y)\leq 6$ can be ignored for large enough  $n, $ and a lower bound for $n$ for this to be true will be given at the end of this proof.

(i)\, Finding $N_{1}(\Sn(T),3)$: \, Let $y$ be in $S_{1}.$  Then 
by Lemma~\ref{lem0001}\, we have $|B_{3}\cap B_{3}(y)|=|B_{2}|+M(y)$ where $M(y)$ is the number of vertices $z\in S_{3}$ with $d(z,y)\leq 3,$  and hence $d(z,y)=2.$ If $z=yt_{1}t_{2}\in S_{3}$ with transpositions $t_{i}$ then also $z^{-1}=t_{2}t_{1}y\in S_{3},$ see Theorem~\ref{con1}. Hence $M(y)$ is the number of all $t_{2}t_{1}y$ belonging to $S_{3}.$ As any
element in $S_{2}$ is of the shape $t_{2}t_{1}$ for some $t_{1}$ and $t_{2}$ we see that $M(y)=|E_{2,3}(y)|=|S_{2}|-|S_{1}|+|S_{0}|$  by (\ref{eij2}). Therefore 
\begin{equation}\label{z0a}
N_{1}(\Sn(T),3)=|B_{2}|+|S_{2}|-|S_{1}|+|S_{0}|=2|S_{2}|+2|S_{0}|.
\end{equation}

(ii)\, Finding $N_{2}(\Sn(T),3)$: \, Let $y=y_{1}y_{2}$ be in $S_{2}$ with transpositions $y_{i}.$ By Lemma~\ref{lem0001}\, we have $|B_{3}\cap B_{3}(y)|=|B_{1}|+M(y)$ where $M(y)$ is the number of vertices $z\in S_{2}\cup S_{3}$ with $d(z,y)\leq 3.$ As above, if $z=yt_{1}\cdots t_{r^{*}}\in S_{2}\cup S_{3}$ with $r^{*}\leq 3,$ consider instead $z^{-1}=t_{r^{*}}\cdots t_{1}y_{2}y_{1}\in S_{2}\cup S_{3},$ that is, all paths from $e_{G}$ of length $\leq 5$ ending in $y_{2}y_{1}$ at a vertex in $S_{2}\cup S_{3}.$ Let $Z$ be the set of all such vertices $z^{-1},$ in particular then  $M(y)=|Z|.$

Let $u=t_{2}t_{1}\in S_{2}$ be arbitrary. If $t_{1}=y_{1}$ then $u=(t_{2}y_{2})y_{2}y_{1}\in Z\cap S_{2}.$ Otherwise $uy_{1}=(t_{2}t_{1}y_{2})y_{2}y_{1}\in Z\cap S_{3}.$ Denote the vertices in $Z$ of this kind by $Z_{0},$ in particular then $|Z_{0}|=|S_{2}|.$

Next let $e=\{v^{+},v^{-}\}\in E_{2,3}(y_{1})$ with $v^{+}=v^{-}y_{1}\in S_{3}$ and $v^{-}=v^{+}y_{1}\in S_{2}.$ Then $v^{-}$ belongs to $Z$ if and only if  $v^{+}y_{2}$ belongs to $S_{2}.$ Let the vertices of this type be denoted by $Z_{1}.$ Thus  $|Z_{1}|$ is the number of $u=v^{+}\in S_{3}$ such that both $uy_{1}$ and $uy_{2}$ belong to $S_{2}.$ When $y_{1}=(1,2)$ and $y_{2}=(2,3)$ then (\ref{e020}) implies that $|Z_{1}|$ is the number of elements in $(1,2,3)(4,5)^{G}$ or $(1,2,3,4)^{G}$ with $1,\,2,\,3$ in the same cycle. This number is 
\begin{equation}\label{z1a}
|Z_{1}|=2{n-3\choose 2}+6(n-3)=(n+2)(n-3)\,\,.\end{equation}   
When $y_{1}=(1,2)$ and $y_{2}=(3,4)$ then $|Z_{1}|$ is the number of elements in $(1,2,3)(4,5)^{G},$ $(1,2,3,4)^{G}$ or $(1,2)(3,4)(5,6)^{G}$ in which $1,\,2$ and $3,\,4$ appear in the same cycle(s). This number is \begin{eqnarray}\label{z1b}|Z_{1}|=4(n-4)+6+{n-4\choose 2}={n\choose 2}\,\,.   \end{eqnarray}

Finally let  $e=\{v^{+},v^{-}\}$ be in $E_{3,4}(y_{1})$ with $v^{+}\in S_{4}$ and $v^{-}\in S_{3}.$ Then $v^{-}$ belongs to  $Z$ if and only if $u=v^{+}\in S_{4}$ has the property  that both $uy_{1}$ and $uy_{2}$ belong to $S_{3}.$ Let $Z_{2}$ be the set of all such vertices $v^{-}$. 
When $y_{1}=(1,2)$ and $y_{2}=(2,3)$ then  $|Z_{2}|$ is the number of elements in $(1,2,3)(4,5,6)^{G},$ $(1,2,3)(4,5)(6,7)^{G},$ $(1,2,3,4)(5,6)^{G}$ or $(1,2,3,4,5)^{G}$ with $1,\,2,\,3$ in the same cycle. This number is
\begin{eqnarray}\label{z2a}
|Z_{2}|&=&4{n-3\choose 3}+{n-3\choose 2}{n-5\choose 2}\nonumber\\
&+&3!(n-3){n-4\choose 2}+4!{n-3\choose 2}\nonumber\\
&=&24{n-3\choose 2}+22{n-3\choose 3}+6{n-3\choose 4}.
\end{eqnarray}
(Note, the term ${n-3\choose 2}{n-5\choose 2}$ accounts for the two choices of a $3$-cycle on $\{1,2,3\}$ while avoiding duplication in the choice of two $2$-cycles from the remaining $n-3$ and $n-5$ vertices, respectively.)

When $y_{1}=(1,2)$ and $y_{2}=(3,4)$ then $|Z_{2}|$ is the number of elements in $(1,2,3)(4,5,6)^{G},$ $(1,2,3)(4,5)(6,7)^{G},$  $(1,2,3,4,5)^{G},$ $(1,2,3,4)(5,6)^{G}$ or $(1,2)(3,4)(5,6)(7,8)^{G}$ in which $1,\,2,$ and $3,\,4$ appear in the same cycle(s). This number is 
\begin{eqnarray}\label{z2b}|Z_{2}|&=&4(n-2)(n-5)+
\left[4(n-4){n-5\choose 2}+  2{n-4\choose 3}\right]   + \nonumber\\
&+&
4!(n-4)+\left[6{n-4\choose 2}+6{n-4\choose 2}+6{n-4\choose 2}\right]+\nonumber\\
&+&\frac12 {n-4\choose2}{n-6\choose 2}\nonumber\\&=&
24(n-4)+(n-5)(13n-44)+
14{n-4\choose 3}+
3{n-4\choose 4}.\,\,\,\,
\end{eqnarray}
  
It is clear that $Z=Z_{0}\cup Z_{1}\cup Z_{2}$ is a disjoint union. Comparing (\ref{z1a})+(\ref{z2a}) to (\ref{z1b})+(\ref{z2b}) one can check that the first expression is bigger than the second for all $n\geq 4.$ Hence $|B_{3}\cap B_{3}(y)|$ takes its maximum when $y\in (1,2,3)^{G}$ for all $n\geq 4.$ Therefore  
\begin{eqnarray}\label{z3a} N_{2}(\Sn(T),\,3)&=&|S_{0}|+|S_{1}|+|S_{2}|+(n+2)(n-3)+\nonumber\\&+&24{n-3\choose 2}+22{n-3\choose 3}+6{n-3\choose 4}\,\,
\end{eqnarray}
and this complete the second part of the theorem. 

We now return to the comment at the beginning of this  proof. Comparing  \,(\ref{z0a})\, to \,(\ref{z3a})\, shows that  $N_{2}(\gG,\,3)> N_{1}(\gG,\,3)$ for all $n\geq 4.$ When $y$ belongs to $S_{3}$ or $S_{4}$ a very rough upper bound for $|B_{3}\cap B_{3}(y)|$ can be obtained by following through the argument in Lemma~\ref{lem00006}.\, By Lemma~\ref{lem0007}\, the downward degree for a vertex in  $S_{j}(e)$  is at most ${j+1\choose 2}$ while the upward degree is at most ${n\choose 2}-j.$ By considering the possible ascent-descent combinations of a path from $y$ to a vertex in $B_{3}$ we can work out that $N_{4}(\gG, 3)\leq -455 + 155{n\choose 2 }$ and 
$N_{3}(\gG, 3)\leq -132 + 65{n\choose 2 }.$ Using the same arguments we can bound $N_{6}(\gG, 3)\leq 1575$ and $N_{5}(\gG, 3)\leq 525.$ Evaluating these inequalities  it can be seen that $N_{2}(\gG,\,3)>N_{j}(\gG,\,3)$ for $j=3,\,4,\,5,\,6$ from $n\geq 16$ onwards. We note that a better lower bound  $n\geq 10$ can be obtained by a more careful albeit tedious count  of the possible paths. This completes the proof.\dne

\medskip
\subsection{ The asymptotic behaviour of $N(\Sn(T),r)$}

The main work in this section will be to find  $N_{1}(\Sn(T),r)$ and $N_{2}(\Sn(T),r)$ for arbitrary $r$ and sufficiently large $n$. It will turn out that this determines   $N(\Sn(T),r)$ in general.
Let $b(n,r)$ denote the cardinality of the ball of radius $r$ in $\Sn(T),$ thus
$$b(n,r)=|B_{r}|=\sum_{0\leq i\leq n}\,\,s_{i}(n)=\sum_{0\leq i\leq n}\,\,c(n,n-i)$$
in terms of the signless Stirling numbers of the first kind. 

First we consider $N_{2}(\Sn(T),r).$  By Lemma~\ref{lem0001} \, we have $N_{2}(\Sn(T),r)=b(n,r-2)+|(S_{r}\cup S_{r-1})\cap B_{r}(y^{*})|$ where   $y^{*}$ suitably is a $3$-cycle or a double transposition. We set $A:=|(S_{r}\cup S_{r-1})\cap B_{r}(y^{*})|$ and  let $y^{*}=y_{1}y_{2}$ with two transpositions $y_{i}.$ We now need to find the number $A$ of all  $z$ in $S_{r}\cup S_{r-1}$ which can be reached on a path $t_{1}t_{2}...t_{r^{*}}$ of length $r^{*}\leq r$ starting from $y^{*}$. This means that $z=y_{1}y_{2}t_{1}t_{2}...t_{r^{*}}$ and applying the inversion automorphism, see \,Theorem~\ref{con1},\, we obtain the element $z^{-1}=t_{r^{*}}\cdots t_{2}t_{1}y_{2}y_{1}$  in $S_{r}\cup S_{r-1}. $ This represents a path from $e=e_{G}$ in which $y_{2}$ and $y_{1}$  are the last edges. If $Z$ denotes the set of  all such elements $z^{-1}$ then $|Z|=A.$

Let $u$ be in $S_{r-1}$ and suppose that $u= t_{r-1}\cdots t_{2}t_{1}$ is the product of  suitable transpositions $t_{i}$. If $t_{1}=y_{1}$ then 
$u= (t_{r-1}\cdots t_{2}y_{2})y_{2}y_{1}$ so that $u\in Z.$ Otherwise $uy_{1}= (t_{r-1}\cdots t_{2}t_{1}y_{2})y_{2}y_{1}$ belongs to $Z.$ Thus every  $u$ in $S_{r-1}$ gives rise to one element in $Z.$ The set of elements of this type is denoted by $Z_{0}$, and in particular $|Z_{0}|=|S_{r-1}|.$ 

The next type of vertices in $Z$  are of  the shape $z^{-1}= t_{r^{*}}\cdots t_{2}t_{1}y_{2}y_{1}$ where both $z^{-1}$ and $t_{r^{*}}\cdots t_{2}t_{1}$ belong to $S_{r-1}$  while $t_{r^{*}}\cdots t_{2}t_{1}y_{2}$ belongs $S_{r}.$ This type will be called $Z_{1},$ evidently this set is disjoint from $Z_{0}.$

The remaining vertices in $Z$ are of the shape $z^{-1}= t_{r^{*}}\cdots t_{2}t_{1}y_{2}y_{1}$ where both $z^{-1}$ and $t_{r^{*}}\cdots t_{2}t_{1}$ belong to $S_{r}$  while $t_{r^{*}}\cdots t_{2}t_{1}y_{2}$ belongs $S_{r+1}.$ These are the vertices of type $Z_{2}$ and it follows that $Z=Z_{0}\cup Z_{1}\cup Z_{2}$ is a disjoint union. Therefore

\begin{equation}\label{z012} 
N_{2}(\Sn(T),r)=b(n,r-1)+\max(|Z_{1}|+|Z_{2}|\,\,:\,\,y^{*}\in S_{2})
\end{equation} 

where  $Z_{1}$ and $Z_{2}$ depend on the choice of $y^{*}$ as either a $3$-cycle or a double transposition. We can now prove the following theorem:

\begin{theorem}
\label{th0071}  \quad Let $\Gamma=\Sn(T)$ be the transposition Cayley graph.  Suppose that $r\geq 2$  and that $y$ is a transposition.  \newline

\vspace{-0.6cm}
(i) \,\,For  $n-1\geq r$  we have 
%\vspace{2.5mm}
\begin{eqnarray}
\label{e025} 
N_{1}(\gG,\,r)&=& b(n,r-1)+ |E_{r-1,r}(y)|=\spa \nonumber\\
&=& 2\cdot\big(|S_{r-1}|+|S_{r-3}|+|S_{r-5}|+\cdots \big)\,\,.\spa%=\nonumber\\\spa&2\,\sum_{0\leq i<\frac r2}\,\,\,c(n,n-r+1+2i)\,\,
\end{eqnarray}
 (ii) For   $n$ sufficiently large we have  
\vspace{1.5mm}
\begin{eqnarray}
\label{e0251} 
N_{2}(\gG,\,r)&>&N_{1}(\gG,\,r)\,\,.%+|S_{r-1}|\,\,.
\end{eqnarray}
\end{theorem}

\pf (i)\, Let $y$ be in $S_{1}$. By Lemma~\ref{lem0001} \,we have $N_{1}(\gG,r)=|B_{r}\cap B_{r}(y)|=
|B_{r-1}|+|S_{r}\cap B_{r}(y)|.$  Hence  we need to find  the number $M(y) =|S_{r}\cap B_{r}(y)|$ of all  $z\in S_{r}$ which can be reached on a path of length $\le r$ from $y$. Such a path is necessarily of the shape $z=yt_{1}t_{2}\cdots t_{r-1},$ consisting of $r-1$ transpositions $t_{i}$. Applying the inversion automorphism, as above, we obtain a path $z^{-1}=t_{r-1}\cdots t_{2}t_{1}y$ starting from $e$ to $z^{-1}\in S_{r}$ with $y$ as last edge. 
Evidently any vertex in $S_{r-1}$ can be reached by a suitable choice of $t_{r-1}\cdots t_{2}t_{1}$ and therefore $M(y)=|E_{r-1,r}(y)|.$ The result now follows from (\ref{eij2}).

(ii)\,  
Let $Z_{1}$ and $Z_{2}$ have the same meaning as in \,(\ref{z012}).   By the first part of this theorem and the equation \,(\ref{z012})\, it will be sufficient to show that $|Z_{2}|\geq |E_{r-1,r}(y)|$ for all large enough $n.$
We evaluate $|Z_{2}|$ for $y^{*}=y_{1}y_{2}$ with  $y_{1}=(1,2)$ and $y_{2}=(2,n).$
The vertices in  $Z_{2}$ are in one-to-one correspondence with vertices $u\in S_{r+1}$ such that $uy_{1}$ and $uy_{2}\in S_{r}.$ By the basic property  \,(\ref{e020})\, this is the case if and only if $1,2,n$ belong to the same cycle of $u$. Let $U$ be the set of such elements, $|U|=|Z_{2}|.$ To count  $|U|$ consider elements $u'$ in the symmetric group $G'$ on $\{1,2,...,n-1\}$ which have the following properties: (i)  both $1$ and $2$ are in the same cycle of $u'$, and (ii) $u'$ has $(n-1) - r$ cycles. Thus $u'$ is in $S_{r}\cap G'$ and it follows from (\ref{eij2}) that the total number of such elements $u'$ is

\vspace{-0.6cm}
\begin{eqnarray}
\label{eij3}a:&=&|\,E_{r-1,r}(y_{1})\,\,\cap\,\, \left\{\,\{x,\,x^{*}\right\}\,\,:\,\,x,\,x^{*}\in G'\}\,|\nonumber\\ &=&|S_{r-1}\cap G'|-|S_{r-2}\cap G'|+|S_{r-3}\cap G'|-...+(-1)^{r-1}|S_{0}\cap G'|\,\,.
%=\nonumber\\ &c\big(n-1,(n-1)-(r-1)\big)-c\big(n-1,(n-1)-(r-1)\big)+....(-1)^{r-1}
\nonumber
\end{eqnarray}

By Lemma~\ref{srn}\, it follows  that the coefficient of the leading power of $n$ in $|S_{r-1}|$ and in $|S_{r-1}\cap G'|$ is the same.
Therefore   

\vspace{-0.6cm}
$$a=|E_{r-1,r}(y_{1})\cap \left\{\{x,x^{*}\}\,:\,x,x^{*}\in G'\right\}|=|S_{r-1}|+f(n)$$ 

for a polynomial $f(n)$ of degree $\leq 2(r-1)-1.$ 
Any $u'$ of the kind just considered is a permutation on $\{1,2,...n-1\}$ in which $1$ and $2$ appear in the same cycle, say of length $c_{u'}\geq 2.$ Now we may insert $n$ into that cycle, in $c_{u'}$ distinct ways, to get  $c_{u'}$ different elements in $U.$ Thus $|Z_{2}|=|U|\geq 2 a$ and therefore   

\vspace{-0.6cm}
\begin{equation} |Z_{2}|\geq 2\cdot|S_{r-1}|+2f(n)\,.\,\end{equation}   

%for  a polynomial $f_{2}(n)$ of degree $\leq 2(r-1)-1.$ 
Since the leading term of $ |E_{r-1,r}(y)|$ is $|S_{r-1}|$ by      \,(\ref{eij2})\,   it follows that  $|Z_{2}|\geq |E_{r-1,r}(y)|$ for all large enough $n.$\dne

\bigskip
In the expression
$N_{2}(\Sn(T),r)=b(n,r-1)+\max(|Z_{1}|+|Z_{2}|\,\,:\,\,y^{*}\in S_{2})$ stated in \,(\ref{z012})\,  the term $|Z_{1}|+|Z_{2}|$ depends on the choice of $y^{*}.$ We therefore turn to evaluating $|Z_{1}|+|Z_{2}|$ for the two possible choices 
 $y^{*}=(1,2,3)$ and  $y^{*}=(1,2)(3,4).$ 

\medskip
This leads us to the following definition. Let    $c_{3^{1}}(n,n-i)$ be the number of vertices in $S_{i}$ in which  the letters $1,2,3$ appear in a single cycle, and let $c_{2^{2}}(n,n-i)$ be the number of vertices in $S_{i}$ in which  the letters $1,2$ and $3,4$ appear in the same cycle or cycles. For instance, 
\begin{eqnarray}c_{3^{1}}(n,n)&=&c_{3^{1}}(n,n-1)=0,\quad c_{3^{1}}(n,n-2)=2,\quad \nonumber \\
c_{3^{1}}(n,n-3)&=&(n+2)(n-3),\,\,\qquad\mbox{and}
\nonumber\\       c_{3^{1}}(n,n-4)&=&24{n-3\choose 2}+22{n-3\choose 3}+6{n-3\choose 4}\,.
\end{eqnarray}
 
 Similarly we have 
 \begin{eqnarray}c_{2^{2}}(n,n)&=& c_{2^{2}}(n,n-1)=0,\quad c_{2^{2}}(n,n-2)=1,\quad \nonumber \\
c_{2^{2}}(n,n-3)&=& {n \choose 2},\,\,\qquad\mbox{and}
\nonumber\\       c_{2^{2}}(n,n-4)&=& 24(n-4)+(n-5)(13n-44)\nonumber\\     &+&
14{n-4\choose 3}+3{n-4\choose 4}.
\end{eqnarray}

For this see again (\ref{z1a}), (\ref{z1b}), (\ref{z2a}) and (\ref{z2b}). As we have already observed in the proof of the last theorem,  the general rule \,(\ref{e020})\, implies that $Z_{1}$ and $Z_{2}$ in \,(\ref{z012})\, satisfy   
\begin{equation}\label{a1a1} |Z_{1}|+|Z_{2}|=c_{3^{1}}(n,n-r) + c_{3^{1}}(n,n-(r+1))\end{equation} 
when  $y^{*}$ is a $3$-cycle and 
\begin{equation}\label{a1a2} |Z_{1}|+|Z_{2}|=c_{2^{2}}(n,n-r) + c_{2^{2}}(n,n-(r+1))\end{equation} when $y^{*}$ is a double transposition.
We obtain an estimate for $c_{3^{1}}(n,n-r)$ and $c_{2^{2}}(n,n-r)$ as follows: 

 \begin{lemma}
\label{lemtwice}  \, For fixed $i\geq 2$ and $n$ sufficiently large   we have 
\begin{equation}\label{twa} c_{3^{1}}(n,n-i)=\frac{2}{(i-2)!}\,{n-3\choose 2}{n-5\choose 2}\cdots{n+3-2i\choose 2} +f_{1}\,\,\,\end{equation}
and 
\begin{equation}\label{twb}c_{2^{2}}(n,n-i)=\frac{1}{(i-2)!}{n-4\choose 2}{n-6\choose 2}\cdots{n+2-2i\choose 2} +f_{2}\,\,\,\end{equation}
where the $f_{i}$ are polynomials in $n$ of degree $<d_{i}=2(i-2).$ 
In particular, $c_{3^{1}}(n,n-i)$ and $c_{2^{2}}(n,n-i)$ are polynomials of degree $d_{i}$ and $c_{3^{1}}(n,n-i)=2\,c_{2^{2}}(n,n-i)+f_{3}$ with a polynomial $f_{3}$ of degree $<d_{i}.$
\end{lemma}

\pf Let $C_{3^{1}}(n,n-i)\subseteq S_{i}$ and $C_{2^{2}}(n,n-i)\subseteq S_{i}$ be the sets counted by $c_{3^{1}}(n,n-i)$ and $c_{2^{2}}(n,n-i)$ respectively. For $i=2,$ when $0!=1,$ we see that $C_{3^{1}}(n,n-2)$ consists of the two $3$-cycles $(1,2,3)$ and $(1,3,2)$ while $C_{2^{2}}(n,n-2)$ consists of the single double transposition $(1,2)(3,4)$ only. This established the base of induction and accounts for the factor $2$ throughout. 
We will prove the statement \,(\ref{twa})\, concerning $c_{3^{1}}(n,n-i),$ the corresponding statement  \,(\ref{twb})\, for $c_{2^{2}}(n,n-i)$ follows in exactly the same way.

If $g\in \Sn$ let ${\rm supp}(g)$ be its {\it support,} that is all symbols moved by $g$. The cardinality of the support of any $g\in C_{3^{1}}(n,n-i)$ is at most $3+2(i-2).$
Let $C_{0}^{i}:=C_{3^{1}}(n,n-i)\cap (1^{n-2i+1}2^{i-2}3^{1})^{G}$ and let $C_{1}^{i}=C_{3^{1}}(n,n-i)\setminus C_{0}^{i}.$ Then $$|C_{0}^{i}|=\frac{2}{(i-2)!}\,{n-3\choose 2}{n-5\choose 2}\cdots{n+3-2i\choose 2} $$ and by induction we assume that $|C_{1}^{i}|=f_{1}$ has degree $<d_{i}$. Since 
$$|C_{0}^{i+1}|=\frac{2}{(i-1)!}\,{n-3\choose 2}{n-5\choose 2}\cdots{n+3-2i\choose 2}\cdot  {n+1-2i\choose 2}$$ it remains to show that  the number of elements 
in $C_{1}^{i+1}$ is a polynomial of degree at most $d_{i}+1.$

%Thus $C_{1}^{i+1}$ is the set of all permutations $g$ in $S_{i+1}$ for which $1,2,3$ are in a cycle of length $\geq 4$ or for which $(1,2,3)$ is a cycle in $g$ while at least one other cycle has length $\geq 3.$
By considering cycle types it is easy to see that any vertex in $C_{1}^{i+1}$ has at least one  neighbour in $C_{0}^{i}$ or in $C_{1}^{i}.$ In the first case, 
if $g=u\cdot(j_{1},j_{2})$ with $u\in C_{0}^{i}$ then at least one of $\{j_{1},\,  j_{2}\}$ must be in the support of $g$ as otherwise $g\in C_{0}^{i+1}.$ The number of such elements $g$ therefore is polynomial  of degree at most $d_{i}+1.$
The number of vertices of the second kind is clearly at most $f_{1}{n\choose 2},$ again of degree at most $d_{i}+1.$ Hence $c_{3^{1}}(n,n-i)$ has the required expression. In the case of $c_{2^{2}}(n,n-i)$ the same arguments apply.\dne

\begin{theorem}
\label{th0009}  \quad Let $\Gamma=\Sn(T)$ be the transposition Cayley graph and suppose that $r\geq 1.$ Then for all sufficiently large $n$  we have 

\vspace{-0.6cm}
\begin{eqnarray}
\label{last} 
N(\gG,r)&=&N_{2}(\gG,r)\,\nonumber\\
&=&b(n,r-1)+c_{3^{1}}(n,n-r) + c_{3^{1}}(n,n-(r+1))\,.\,\,\,\spa
\end{eqnarray}
\end{theorem}

{\sc Remark:}\quad For $r\leq 3$ we already have computed the value of $N(\Sn(T),r)$ in Theorems~\ref{lem0008},\,\,\ref{th005}\,\, and \ref{th006}\,\, when $N(\Sn(T),r)$ indeed agrees with (\ref{last}). In these theorems the lower bound on $n$ was explicit and hence better than the condition here. Nevertheless, analysing the arguments here it is likely that the bound $n>3r$ is sufficient. 

\pf We can assume that $r>3$. By \,(\ref{z012}),\, (\ref{a1a1})\, and   Lemma~\ref{lemtwice} \, it follow that $N_{2}(\gG,r)=b(n,r-1)+c_{3^{1}}(n,n-r) + c_{3^{1}}(n,n-(r+1))$ and this establishes  the second equation.  By Theorem~\ref{th0071} \,we know that $N_{2}(\Sn(T),r)>N_{1}(\Sn(T),r)$ and both terms are polynomial of degree $2(r-1).$ It remains to show that $N_{s}(\Sn(T),r)$ is polynomial of degree $<2(r-1)$ for $2<s.$ For $y\in S_{s}$ consider all paths $z=yt_{1}\cdots t_{r^{*}}$ of length $r^{*}\leq r$ to a vertex $z\in S_{s^{*}}$ with $s^{*}\leq r.$ If $a$ and $b$ are the number of ascents and descents then $a+b=r^{*}\leq r$ and $a-b=s^{*}-s.$ From this it follows that $a<r-1$ and the required fact now follows from Lemma~\ref{lem00006}.\dne

\end{document}